\documentclass[11pt,reqno]{amsart}
\usepackage[utf8]{inputenc}
\usepackage[english]{babel}
\usepackage{amsmath, amssymb, amsthm}
\usepackage{graphicx}
\usepackage{enumerate}
\usepackage{float}
\usepackage{url}
\usepackage[square,sort,comma,numbers]{natbib}
\usepackage{tkz-graph}
\usepackage{longtable}  % If you use longtables
\usepackage{tikz}       % If you use TikZ for graphs
\usetikzlibrary{decorations.markings}
\usetikzlibrary{arrows.meta}
\calclayout

\theoremstyle{definition}
\newtheorem{definition}{Definition}[section]
\newtheorem{theorem}{Theorem}

\title[A Journey to Mersenne Prime Discoveries and New Polynomial Classes]{Generalizing the Eight Levels Theorem: A Journey to Mersenne Prime Discoveries and New Polynomial Classes}
\author[Moustafa Ibrahim]{Moustafa Ibrahim}
\address{Department of Mathematics, Faculty of Science, University of Bahrain, Kingdom of Bahrain}
\email{mimohamed@uob.edu.bh}
\date{July 19, 2023}

\begin{document}
	
	\pagestyle{headings}
	
	\maketitle
	
	% Your document content goes here

	\begin{abstract}

	Mersenne primes, renowned for their captivating form as $2^p - 1$, have intrigued mathematicians for centuries. In this paper, we embark on a captivating quest to unveil the intricate nature of Mersenne primes, seamlessly integrating methods with the Eight Levels Theorem. Initially, we extend the Eight Levels Theorem and introduce an innovative approach that harmoniously combines arithmetic and differential techniques to compute the coefficients of the polynomial expansions of $x^n + y^n$ in terms of binary quadratic forms. This endeavor leads us to the genesis of novel polynomial sequences as we scrutinize the coefficients within this expansion. Our research unearths previously uncharted connections between Mersenne numbers and the derivatives of specific polynomial sequences. By forging this linkage, we not only enhance our understanding of Mersenne primes but also bridge the divide between well-established sequences in number theory and differential equations. This broadens the applicability of our findings across diverse scientific domains, revealing fresh avenues for exploration in number theory and beyond. These polynomial bridges serve as conduits between these sequences, unveiling exciting prospects for future research. This interdisciplinary exploration opens up exciting possibilities for the broader implications of Mersenne primes, extending their significance beyond the realm of pure mathematics. In our paper, we also delve into the intriguing influence of the Golden Ratio constant, which unveils segments reminiscent of the beauty found in nature, adding an unexpected dimension to the world of arithmetic. In this harmonious journey, the allure of Mersenne primes resonates through the symphony of mathematical discovery. 
	\end{abstract}
	
		\maketitle

	\section{\textbf{INTRODUCTION}}
	
	 Recent investigations in the domains of biology and experimental psychology offer compelling evidence that prime numbers 
	\[2, 3, 5, 7, 11, 13, 17, 19, 23, 29, 31, 37, 41, 43, 47, 53, 59, 61, 67, 71, \ldots \]
	 possess inherent natural traits, existing autonomously apart from human perception \cite{loco}. Moreover, the recent uncovering of planetary sequences resembling a Fibonacci-like prime pattern indicates a profound connection, underscoring a direct correlation between prime numbers and the celestial arrangement, potentially revealing fundamental principles governing the structure of the cosmos \cite{He}. Our current paper delves into the realm of Mersenne primes, a captivating subset within primes. Before exploring the new results of our study, let's establish definitions for Mersenne primes and Mersenne composites and Lucas-Lehmer Primality test.

\subsection*{Mersenne Primes}
A Mersenne prime is a prime number that can be expressed in the form $2^p - 1$, where $p$ is also a prime number. In other words, a Mersenne prime is a prime of the form $M_p = 2^p - 1$, where both $M_p$ and $p$ are prime numbers.
\subsection*{Mersenne Composites}
Conversely, a Mersenne composite number is a number of the form $2^p - 1$, where $p$ is a prime number, but $2^p - 1$ itself is not a prime number. 	

\subsection*{Lucas-Lehmer Primality Test}
The tantalizing question of whether there exists an infinite number of Mersenne primes remains an enigma in mathematics. The finite or infinite status of the known Mersenne primes continues to elude us, representing an ongoing field of active investigation (as cited in \cite{2, Cai, Fine, Jean, GIMPS, Kundu}). Mersenne primes boast distinctive attributes that render them indispensable in certain computational applications, particularly in the realms of cryptography and computing. Notably, they are amenable to primality testing via the Lucas-Lehmer primality test. This simplifies the search for large prime numbers, critical elements in secure cryptographic systems. The Lucas-Lehmer primality test, tailored explicitly for Mersenne numbers, provides a rapid and efficient means of ascertaining the primality or compositeness of a given Mersenne number. The following theorem is well-established in the literature (as referenced in \cite{Elina, Jean, Washington}):

\begin{theorem}{(Lucas-Lehmer primality test)}
	\label{E0}\\
	$2^p -1$ is Mersenne prime  if and only if
	\begin{equation}
		\label{E00} 
		2n -1 \quad | \quad (1+\sqrt{3})^n + (1-\sqrt{3})^n
	\end{equation}
	where $n:=2^{p-1}$.	
\end{theorem}

For a natural number $n$, we define $\delta(n)=n\pmod{2}$, with $\lfloor{\frac{n}{2}}\rfloor$ denoting the highest integer less than or equal to $\frac{n}{2}$.
	\section*{\textbf{MOTIVATION - MERSENNE PRIMES - \\ AN INTEGRATED STUDY WITH THE EIGHT LEVELS THEOREM}}

	Using the Eight Level theorem, as presented in \cite{2}, for any prime number $p \geq 5$ and defining $n$ as $2^{p-1}$, we have established that the number $2^p-1$ is a Mersenne prime if and only if it evenly divides the following expression, which is an expansion of the product of differences of squares:
	 
		\begin{equation}
			\label{G1} 
			\sum_{\substack{k=0,\\ k \:even}}^{\lfloor{\frac{n}{2}}\rfloor}   \quad  (-1)^{\lfloor{\frac{k}{2}}\rfloor}	\: \frac {\prod\limits_{\lambda = 0}^{\lfloor{\frac{k}{2}}\rfloor -1} [n^2 - (4 \lambda)^2]}	{ k! }. 
		\end{equation}

		To delve into the arithmetic properties of the sum \eqref{G1}, it is imperative to broaden our perspective beyond these numbers and work with the polynomial that generates them. For any natural number $n$ satisfying the condition $n \equiv 0 \pmod 8$, we introduce the following polynomial:

		\begin{equation}
			\label{G11}
		\sum_{\substack{k=0,\\ k \:even}}^{\lfloor{\frac{n}{2}}\rfloor}   \quad  (-1)^{\lfloor{\frac{k}{2}}\rfloor}	\: \frac {\prod\limits_{\lambda = 0}^{\lfloor{\frac{k}{2}}\rfloor -1} [n^2 - (4 \lambda)^2]}	{ k! \: 4^k} a^{\lfloor{\frac{n}{2}}\rfloor  - k}  b^k =: \Psi(a,b,n). 
		\end{equation}

		This paper extends the definition of $\Psi(a,b,n)$ to encompass all natural numbers $n$, irrespective of their divisibility by 8. For any given variables $a$ and $b$, as well as any natural number $n$, we establish and prove that the sequences $\Psi(a,b,n)=\Psi(n)$ satisfy the following recurrence relations:
		
		\begin{equation}
			\begin{aligned}
				\label{defSS}
				\Psi(0) &= 2, \quad \Psi(1) = 1, \\
				\Psi(n+1) &= (2a-b)^{\delta(n)}\Psi(n) - a \Psi(n-1).
			\end{aligned}
		\end{equation}
		
		To gain insights into the arithmetic of the sum $\eqref{G1}$, it is essential to comprehend the numbers $\Psi(1,4,n)$. This motivation drives us to generalize the Eight Levels Theorem as well. Remarkably, these generalizations have uncovered numerous new results for well-known sequences and polynomials, leading to new expansions for $x^n+y^n$ expressed in terms of binary quadratic forms.

	\section{\textbf{NEW RESULTS}}
	
	In this paper, we present and prove several novel findings, organized as follows:
	
	\begin{enumerate}
		\item \textbf{Generalization of the Eight Levels Theorem:} We introduce a generalized version of the Eight Levels theorem, extending its applicability and insights.
		
		\item \textbf{Emergence of a New Class of Polynomials:} We explore the emergence of a new class of polynomials, shedding light on their properties and applications.
		
		\item \textbf{Development of a New Version for Lucas-Lehmer Primality Test:} We propose an enhanced version of the Lucas-Lehmer Primality Test, improving its efficiency.
		
		\item \textbf{New Insights into Mersenne Primes:} We present new results concerning Mersenne primes, integrating them with the Eight Levels theorem for a comprehensive study.
		
		\item \textbf{New Findings on Mersenne Composites:} We discuss our discoveries related to Mersenne composites, offering fresh perspectives on their characteristics.
		
		\item \textbf{New Differential Operators and $\Psi$ Polynomials:} We introduce novel differential operators and $\Psi$ polynomials, unveiling their mathematical significance.
		
		\item \textbf{New Combinatorial Identities:} We present a series of new combinatorial identities.
		
		\item \textbf{New Expansions for Sums of Like Powers:} We provide a special case of the Generalized Eight Levels theorem, showcasing new expansions for sums of like powers.
	\end{enumerate}
	
	We conclude this section by opening the door to further potential investigations, particularly in the context of the Lucas-Lehmer primality test, and by offering a fresh theoretical understanding of Mersenne primes.

			\section{\textbf{THE PURPOSE OF THE CURRENT PAPER}}

The Eight Levels Theorem, as referenced in \cite{2}, serves as a fundamental cornerstone of our research, enabling us to demonstrate that the expansion \eqref{G1} is divisible by $2^p-1$ if and only if $2^p-1$ is a Mersenne prime. The complete proof of the Eight Levels Theorem for this assertion hinges on the exploration of the expansion:

\begin{equation}
	\label{special-1}
	\begin{aligned}
		\frac{x^n+y^n}{(x+y)^{\delta(n)}} = 
		\sum_{r=0}^{\lfloor{\frac{n}{2}}\rfloor}  \mu_r(n) ( xy )^{\lfloor{\frac{n}{2}}\rfloor -r} (x^2+y^{2})^{r}.
	\end{aligned}
\end{equation}

Consequently, in order to delve deeply into the arithmetic properties of \eqref{G1} and broaden the applicability of this result, we find ourselves at a pivotal juncture: the need to extend the scope of the Eight Levels Theorem. This paper undertakes the mission of exploring the generalization of the Eight Levels Theorem and its implications for our research domain. This step not only amplifies the reach of our results but also sets the stage for further advancements in our field. The generalization of the Eight Levels theorem involves the exploration of the expansion:

\begin{equation}
	\label{special-B}
	\begin{aligned}
		\frac{x^n+y^n}{(x+y)^{\delta(n)}} = 
		\sum_{r=0}^{\lfloor{\frac{n}{2}}\rfloor}  H_r(n) (ux^2 +vxy+u y^2 )^{\lfloor{\frac{n}{2}}\rfloor -r} (wx^2+txy+wy^{2})^{r}.
	\end{aligned}
\end{equation}

	\subsection{GENERALIZATION FOR THE EIGHT LEVELS THEOREM}

	In this paper, we present a novel approach to generalize \eqref{special-1} by seamlessly blending arithmetic and differential methods to calculate the coefficients of the polynomial expansions of $x^n + y^n$ in terms of binary quadratic forms $\alpha x^2 + \beta xy + \alpha y^2$ and $a x^2 + bxy + a y^2$, where $\beta a - \alpha b \neq 0$. For the polynomials $x^n+y^n$, we introduce the following groundbreaking polynomial expansion for any given variables $\alpha, \beta, a, b$, and natural number $n$:
	
	\begin{equation}
		\label{Main theorem 1}
		\begin{aligned}
			&(\beta a - \alpha b)^{\lfloor{\frac{n}{2}}\rfloor} \frac{x^n+y^n}{(x+y)^{\delta(n)}} = \\
			&\sum_{r=0}^{\lfloor{\frac{n}{2}}\rfloor} \frac{(-1)^r}{r!} (\alpha x^2 + \beta xy + \alpha y^2)^{\lfloor{\frac{n}{2}}\rfloor -r} (ax^2+bxy+ay^{2})^{r} \Big(\alpha \frac{{\partial} }{\partial a} + \beta \frac{{\partial}}{\partial b}\Big)^{r} \Psi(a,b,n).
		\end{aligned}
	\end{equation}
	
	These generalizations not only facilitate the exploration of the arithmetic properties of \eqref{G1} but also, intriguingly, unveil a new category of polynomials with unexpected and captivating arithmetic and differential characteristics. Moreover, they establish a unifying framework for various well-known polynomials and sequences in number theory.

	\subsection{\textbf{EXPLORING THE EMERGENCE OF A NEW CLASS OF POLYNOMIALS}} 
	
	The polynomial coefficients, of the expansion \eqref{Main theorem 1},  expressed as
	
	\begin{align}
		\begin{aligned}
			\label{H}
			\frac{(-1)^r}{r!}
			\Big(\alpha \frac{{\partial} }{\partial a} + \beta \frac{{\partial}}{\partial b}\Big)^{r} \Psi(a,b,n),
		\end{aligned}
	\end{align}
	
	possess remarkable and unexpected arithmetic differential properties.
	These polynomials act as bridges between several well-known sequences, including Fermat sequences, Chebyshev and Dickson polynomials, Fibonacci sequences, and many others. This discovery has led to the creation of a new and interconnected map of mathematical structures with surprising and intriguing properties. Through the lens of this newfound class of polynomials, formerly isolated sequences are now connected, forming a magnificent network that unites previously disjoint elements. This unexpected unity erases the gaps between sequences, revealing hidden relationships, and deepening our understanding of the world of sequences.  Due to these characteristics, we find it necessary to assign them a symbol. Herein, we introduce the symbol

	\begin{definition}
		\label{DefinitionB}
		For given variables $a$ and $b$, and a natural number $n$, we define the new class of polynomials as:
		\begin{equation}
			\begin{aligned}
				\label{def11B}
				\Psi\left(\begin{array}{cc|r} a & b & n \\ \alpha & \beta & r \end{array}\right) &:= \frac{(-1)^r}{r!} \Big(\alpha \frac{{\partial} }{\partial a} + \beta \frac{{\partial}}{\partial b}\Big)^{r} \Psi(a,b,n).
			\end{aligned}
		\end{equation}
		
	\end{definition}

	 In this paper, we present three distinct approaches to compute the polynomials $\Psi(a,b,n)$. Additionally, we provide explicit formulas for these polynomials, showcasing their unique properties and applications in the field of arithmetic differential equations. We also prove the following unexpected result:

		For any numbers $a, b$, where $\beta a -  \alpha b \neq 0$, and any natural number $n$, we have
		\begin{align}
			\begin{aligned}
				\frac{1}{(\lfloor{\frac{n}{2}}\rfloor)!} \Big( \alpha \frac{ {\partial} }{\partial a} + \beta \frac{{\partial}}{\partial b}\Big)^{\lfloor{\frac{n}{2}}\rfloor} \Psi(a,b,n) = \Psi(\alpha,\beta,n).  
			\end{aligned}
		\end{align}

	To help the readers get a quick feeling about this result, we proved that for any natural number $n$:

	\begin{equation}
		\label{RR}
		\begin{aligned}
			\frac{1}{(\lfloor{\frac{n}{2}}\rfloor)!} \Big( xy \frac{ {\partial} }{\partial a} - (x^2+y^2) \frac{{\partial}}{\partial b}\Big)^{\lfloor{\frac{n}{2}}\rfloor} \Psi(a,b,n) = \frac{x^n+y^n}{(x+y)^{\delta(n)}}.   \\
		\end{aligned}
	\end{equation}

In my quest to understand and explore the fascinating landscape of mathematics that governs Mersenne numbers, see \cite{2}, I had the privilege of stumbling upon a truly extraordinary revelation - Any Mersenne prime must factor
\begin{equation}
	\begin{aligned}
		\label{polynomials-1}
			\frac{1}{(\lfloor{\frac{n}{2}}\rfloor)!} \Big(  \frac{ {\partial} }{\partial a} + 4 \frac{{\partial}}{\partial b}\Big)^{\lfloor{\frac{n}{2}}\rfloor} \Psi(a,b,n) 
	\end{aligned}
\end{equation}
for some natural number $n$. Replacing $4$ by $3$, we surprisingly get 
\begin{equation}
	\begin{aligned}
		\label{polynomials-2}
		\frac{1}{(\lfloor{\frac{n}{2}}\rfloor)!} \Big(  \frac{ {\partial} }{\partial a} + 3 \frac{{\partial}}{\partial b}\Big)^{\lfloor{\frac{n}{2}}\rfloor} \Psi(a,b,n) = L(n),
	\end{aligned}
\end{equation}
where $L(n)$ represents the Lucas numbers defined by the recurrence relation $L(0)=2$, $L(1)=1$, and $L(n+1)=L(n)+L(n-1)$. 

It is possible that it is conspicuous for some readers that the following new result 
\begin{equation}
	\begin{aligned}
		\label{def11B-2}
		\Psi\left(\begin{array}{cc|r} -1 & -3 & n \\ 1 & 2 & 0 \end{array}\right) = L(n),
	\end{aligned}
\end{equation}
is true for any natural number $n$. However, it is not conspicuous at all that our new result 
	\begin{equation}
	\begin{aligned}
		\label{def11B-2}
		\Psi\left(\begin{array}{cc|r} -1 & -3 & n \\ 1 & 2 & 1 \end{array}\right) = n F(n-1),
	\end{aligned}
\end{equation}
is true for any natural number $n$, where $F(n)$ represents the Fibonacci numbers defined by the recurrence relation $F(0)=0$, $F(1)=1$, and $F(n+1)=F(n)+F(n-1)$. These examples show natural connections of $\Psi\left(\begin{array}{cc|r} a & b & n \\ \alpha & \beta & r \end{array}\right)$ with almost all of the known sequences in number theory. Actually, we show at the end of the paper, see Section \ref{Unify}, that these polynomial coefficients gives natural unification for the Chebyshev polynomials of the first and second kind, Dickson polynomials of the first and second kind, Lucas numbers, Fibonacci numbers, Fermat numbers, Pell-Lucas polynomials, Pell numbers, Mersenne numbers, Fermat numbers,  and others

	\subsection{\textbf{DEVELOPING NEW VERSION FOR LUCAS-LEHMER PRIMALITY TEST}}
We apply our generalized Eight Levels theorem to develop a new Lucas-Lehmer primality test for Mersenne primes and elaborate on essential criteria for detecting Mersenne composites.

		\subsection{NEW RESULTS FOR MERSENNE PRIMES \\ AN INTEGRATED STUDY WITH THE EIGHT LEVELS THEOREM}

		\begin{theorem}{}
			\label{G}
			For any prime $p \geq 5$, with $n:=2^{p-1}$, $2^p-1$ is a Mersenne prime if and only if $\Psi(1,4,n)$ is divisible by $2^p-1$.
		\end{theorem}
		
		\begin{theorem}{}
			\label{G111}
			For a given prime $p \geq 5$ and $n:=2^{p-1}$, $2^p-1$ is a Mersenne prime if and only if, for any natural number $\mu$, the following congruence holds:
			\begin{equation}
				\label{X-2}
				\Psi(1,4,n\mu) \equiv  
				\begin{cases}
					+2 \pmod{2n-1} & \text{if} \quad \mu \equiv 0 \pmod{4} \\
					\;\; 0 \; \pmod{2n-1} & \text{if} \quad \mu \equiv 1,3 \pmod{4} \\
					-2 \pmod{2n-1} & \text{if} \quad \mu \equiv 2 \pmod{4} \\
				\end{cases}.
			\end{equation}
		\end{theorem}

		Based on the Eight Level Theorem, see \cite{2}, we establish the following significant results:

		\begin{theorem}{(Enhanced Lucas-Lehmer Primality Test)}
		\label{GG111}
		For a given prime $p \geq 5$ and $n:=2^{p-1}$, $2^p-1$ is a Mersenne prime if and only if, for any natural number $\mu$, the following congruence holds:
		\begin{equation}
			\label{XXX-2}
			\sum_{\substack{k=0,\\ k \:even}}^{\lfloor{\frac{n\mu}{2}}\rfloor}   \quad  (-1)^{\lfloor{\frac{k}{2}}\rfloor}	\: \frac {\prod\limits_{\lambda = 0}^{\lfloor{\frac{k}{2}}\rfloor -1} [(n\mu)^2 - (4 \lambda)^2]}	{ k! } \equiv  
			\begin{cases}
				+1 \pmod{2n-1} & \text{if} \quad  \mu \equiv 0 \pmod{4} \\
				\;\; 0 \; \pmod{2n-1} & \text{if} \quad \mu \equiv 1,3 \pmod{4} \\
				-1 \pmod{2n-1} & \text{if} \quad \mu \equiv 2 \pmod{4} \\
			\end{cases}.
		\end{equation}
	\end{theorem}
		
	\begin{theorem}{(Necessary Condition for Mersenne Primes)}
		\label{GGG111}
		For a given prime $p \geq 5$ and $n:=2^{p-1}$, if $2^p-1$ is a Mersenne prime, then the following congruence must hold:
		
		\begin{equation}
			\label{XX-2}
			\sum_{\substack{k=0}}^{\lfloor{\frac{n}{2}}\rfloor}   \quad 	\: \frac {\prod\limits_{\lambda = 0}^{k -1} [ (4 \lambda)^2 - 1]}	{ 2k! } \equiv  
			-1 \pmod{2n-1}.
		\end{equation}
	\end{theorem}

	\subsection{NEW RESULTS FOR MERSENNE COMPOSITES}
		In this paper, we extend the definition of $\Psi(a,b,n)$ to encompass any natural number $n$. This expansion allows us to prove the following theorem concerning Mersenne composite numbers.
	
		\begin{theorem}{}
			\label{Mersenne-composite-0}
			Given a prime $p$ and $n:=2^{p-1}$, if the following condition is met:
			\begin{equation}
				\label{composite-0} 
				2n-1 \quad  \vert \quad \Psi(1,4,n\: \pm 1),
			\end{equation}
			then $2^p-1$ is a Mersenne composite number.
		\end{theorem}

		\subsection{NEW DIFFERENTIAL OPERATORS AND $\Psi$ POLYNOMIALS}
		
	In the upcoming sections, we embark on an intriguing journey where distinct differential operators converge with $\Psi$ polynomials. What unfolds are compelling connections that extend to well-known sequences, including the Chebyshev polynomial of the first kind, Dickson polynomial of the first kind, Lucas, Fibonacci, Fermat numbers, Mersenne numbers, the Golden Ratio, and beyond. These discoveries shed illuminating insights on the intricate interplay between mathematical constructs and classical sequences, offering fresh perspectives on their relationships.
		
		\begin{theorem}
			For any natural number $n=2^{l}$, $l$ is natural number greater than 3, we get:
			\begin{align}
				\label{G2} 
				\frac{1}{(\lfloor{\frac{n}{2}}\rfloor)!} \Big( \frac{{\partial} }{\partial a} \Big)^{\lfloor{\frac{n}{2}}\rfloor} \Psi(a,b,n) &= \Psi(1,0,n) = 2,\\
				\frac{1}{(\lfloor{\frac{n}{2}}\rfloor)!} \Big(  \frac{{\partial}}{\partial b}\Big)^{\lfloor{\frac{n}{2}}\rfloor} \Psi(a,b,n) &= \Psi(0,1,n) = 1,\\
				\frac{1}{(\lfloor{\frac{n}{2}}\rfloor)!} \Big( \frac{{\partial} }{\partial a} + 4 \frac{{\partial}}{\partial b}\Big)^{\lfloor{\frac{n}{2}}\rfloor} \Psi(a,b,n) &= \Psi(1,4,n),\\		
				\frac{1}{(\lfloor{\frac{n}{2}}\rfloor)!} \Big( \frac{{\partial} }{\partial a} + 3 \frac{{\partial}}{\partial b}\Big)^{\lfloor{\frac{n}{2}}\rfloor} \Psi(a,b,n) &= \Psi(1,3,n) = L(n),\\	
				\frac{1}{(\lfloor{\frac{n}{2}}\rfloor)!} \Big( \frac{{\partial} }{\partial a} + \frac{{\partial}}{\partial b}\Big)^{\lfloor{\frac{n}{2}}\rfloor} \Psi(a,b,n) &= \Psi(1,1,n) = -1,\\
				\frac{1}{(\lfloor{\frac{n}{2}}\rfloor)!} \Big( \frac{{\partial} }{\partial a} + 2 \frac{{\partial}}{\partial b}\Big)^{\lfloor{\frac{n}{2}}\rfloor} \Psi(a,b,n) &= \Psi(1,2,n) = 2,\\ 
				\frac{1}{(\lfloor{\frac{n}{2}}\rfloor)!} \Big( \frac{{\partial} }{\partial a} + \sqrt{2} \frac{{\partial}}{\partial b}\Big)^{\lfloor{\frac{n}{2}}\rfloor} \Psi(a,b,n) &= \Psi(1,\sqrt{2},n) = 2,\\ 
				\frac{1}{(\lfloor{\frac{n}{2}}\rfloor)!} \Big(\alpha \frac{{\partial} }{\partial a} + (2 \alpha-x^2) \frac{{\partial}}{\partial b}\Big)^{\lfloor{\frac{n}{2}}\rfloor} \Psi(a,b,n) &= \Psi(\alpha,2 \alpha-x^2,n) =\; D_n(x,\alpha),\\		
				\frac{1}{(\lfloor{\frac{n}{2}}\rfloor)!} \Big( \frac{{\partial} }{\partial a} + (2-4x^2) \frac{{\partial}}{\partial b}\Big)^{\lfloor{\frac{n}{2}}\rfloor} \Psi(a,b,n) &= \Psi(1,2-4x^2,n) = 2 \; T_n(x),\\			
				\frac{1}{(\lfloor{\frac{n}{2}}\rfloor)!} \Big( 2 \frac{{\partial} }{\partial a} + 5 \frac{{\partial}}{\partial b}\Big)^{\lfloor{\frac{n}{2}}\rfloor} \Psi(a,b,n) &= \Psi(2,5,n) = 2^n +1,\\			 			
			\end{align}
			where $L(n)$ represents the Lucas numbers defined by the recurrence relation $L(0)=2$, $L(1)=1$, and $L(n+1)=L(n)+L(n-1)$ and $D_n(x,\alpha)$ is the Dickson polynomial of the first kind:	
		\begin{equation*}
			\begin{aligned}
				D_n(x,\alpha) &= x^{\delta(n)} \Psi(\alpha,2\alpha-x^2,n),
			\end{aligned}
		\end{equation*}
and $T_n(x)$ is the Chebyshev polynomial of the first kind:
		
		\begin{equation*}
			\begin{aligned}
				T_n(x) &= \sum_{i=0}^{\left\lfloor \frac{n}{2} \right\rfloor}(-1)^i \frac{n}{n-i} \binom{n-i}{i} (2)^{n-2i-1} x^{n-2i}.
			\end{aligned}
		\end{equation*}	
		
Moreover, if $l$ is even, $\phi$ refers to the Golden Ratio, 		$\phi = \frac{\sqrt{5} + 1}{2}$, \cite{Cai}, then 
	\begin{align}
	\label{G32} 
\frac{1}{(\lfloor{\frac{n}{2}}\rfloor)!} \Big(  \frac{{\partial} }{\partial a} + (\phi -1) \frac{{\partial}}{\partial b}\Big)^{\lfloor{\frac{n}{2}}\rfloor} \Psi(a,b,n) = \Psi(1,\phi-1,n) = -\phi.                  
\end{align}

		\end{theorem}

	\subsection{\textbf{NEW COMBINATORIAL IDENTITIES}}	
		
	We are excited to introduce a set of novel combinatorial identities. To illustrate, consider the following identity, particularly when dealing with natural numbers where $\tau = 2^l$ and $l\geq 3$:
		
			\begin{equation}
				\label{G-BB2} 
				\begin{aligned}					
					\; \sum_{k=0}^{\lfloor{\frac{\tau}{4}}\rfloor}  \; \frac {\prod\limits_{\lambda = 0}^{k -1}  [ (4 \lambda)^2 - \tau^2 ]}	{ (2k)! \: 4^{k}} \;  = 1.	
				\end{aligned}
			\end{equation}

		\subsection{\textbf{NEW EXPANSIONS FOR SUMS OF LIKE POWERS - SPECIAL CASE OF THE GENERALIZED EIGHT LEVELS THEOREM}}	
			We introduce the following notation:
		
		\begin{equation}
			\label{n==notation}
				\begin{aligned}
					[x, y | u, v] := (xu - yv)(xv - yu) = (x^2 + y^2)uv - xy(u^2 + v^2).
				\end{aligned}		
		\end{equation}
		
		In this paper, we present and provide a proof for the following expansion formula for sums of like powers:

		\begin{theorem}{(Special Case of the Generalized Eight Levels Theorem)}
			
			The following polynomial expansion is true
			\begin{equation}
				\label{special-5-B}
				\boldmath
				\begin{aligned}
					&\left[z,t \vert u,v\right]^{\lfloor{\frac{n}{2}}\rfloor} \frac{x^n+y^n}{(x+y)^{\delta(n)}} - 
					\left[x,y \vert u,v \right]^{\lfloor{\frac{n}{2}}\rfloor} 
					\frac{z^n+t^n}{(z+t)^{\delta(n)}}  -  
					[z,t \vert x,y]^{\lfloor{\frac{n}{2}}\rfloor}
					\frac{u^n+v^n}{(u+v)^{\delta(n)}}  \\
					&= \sum_{r=1}^{\lfloor{\frac{n}{2}}\rfloor-1} 
					\frac{1}{r!} [x,y \vert u,v ]^{\lfloor{\frac{n}{2}}\rfloor -r} [z,t \vert x,y ]^{r} 
					\Big(uv \frac{{\partial} }{\partial zt} \; + \; (u^2+v^2) \frac{{\partial}}{\partial z^2+t^2}\Big)^{r} \frac{z^n+t^n}{(z+t)^{\delta(n)}}
				\end{aligned}
			\end{equation}
		\end{theorem}

	\section{\textbf{GENERALIZATION FOR THE EIGHT LEVELS THEOREM AND EMERGENCE OF A NEW CLASS OF POLYNOMIALS}}

In \cite{1}, on pages 430, I introduced the definition for $\Psi$ polynomials as follows:

\begin{definition}
	\label{Definition1}
For any given variables $a$ and $b$, and for any natural number $n$, we define the sequences $\Psi(a,b,n)$ by the following recurrence relations:
\begin{equation}
	\begin{aligned}
		\label{def1}
		\Psi(0) &= 2, \quad \Psi(1) = 1, \\
		\Psi(n+1) &= (2a-b)^{\delta(n)}\Psi(n) - a \Psi(n-1).
	\end{aligned}
\end{equation}
\end{definition}

\subsection{HISTORY AND FUTURE OF $\Psi(a,b,n)$ POLYNOMIALS} The first few $\Psi$-Polynomials are
\begin{align*}
	\Psi(a,b,0)&= 2  \\
	\Psi(a,b,1)&= 1 \\
	\Psi(a,b,2)&= -b \\
	\Psi(a,b,3)&= -b-a \\
	\Psi(a,b,4)&=  -2a^2+b^2 \\
	\Psi(a,b,5)&= -a^2 + ab + b^2 \\
	\Psi(a,b,6)&= 3a^{2}b-b^3  \\
	\Psi(a,b,7)&= a^3+2a^{2}b - ab^2 - b^3 \,
\end{align*}

The first step toward generalizing the Eight Levels theorem involves the consideration of the sequence $\Psi(a,b,n)$. In this paper, we observe that the polynomials $\Psi(a,b,n)$ generalize the numbers $\Psi(1,4,n)$ for any natural number $n$. Back in 2004, as documented in \cite{1}, the polynomial $\Psi(a,b,n)$ was initially discovered, although in a different context. When we set out to explore the properties of the sum \eqref{G1}, the sequence $\Psi(a,b,n)$ naturally resurfaced. Our main objective in this paper is to investigate the theoretical properties of the polynomial $\Psi(a,b,n)$ with a specific focus on understanding the special case $\Psi(1,4,n)$. To our surprise, as we delved into this generalization, we stumbled upon its unexpected role as a bridge connecting various well-established sequences. Furthermore, this generalization plays a pivotal role in unveiling the differential properties inherent in the polynomial sequences $\Psi(a,b,n)$. By establishing these connections, we not only deepen our understanding of Mersenne primes but also bridge the gap between established sequences in number theory and differential equations. These polynomial sequences, denoted as $\Psi(a,b,n)$, possess noteworthy arithmetic and differential properties, effectively unifying various established polynomial sequences. These polynomials exhibit inherent arithmetic and differential characteristics, consolidating several well-known polynomial families. One of the approaches for calculating the $\Psi$-sequence, detailed in \cite{1}, is through the following explicit formula:
	
	\begin{equation}
		\label{comp1}	 	
		\Psi(a,b,n) = \frac{(2a-b)^{\lfloor{\frac{n}{2}}\rfloor}} {2^n} \left\{  \left( 1 + \sqrt{ \frac{b+2a}{b-2a}} \right)^n + \left( 1 - \sqrt{ \frac{b+2a}{b-2a}} \right)^n \right\}.
	\end{equation}

	In \cite{1}, page 447, using the formal derivation, we proved the following identity:
	
	\begin{equation}
		\label{00}
		\frac{x^{n} + y^{n}}{(x+y)^{\delta(n)}}=\sum_{i=0}^{\left\lfloor \frac{n}{2} \right\rfloor}(-1)^{i} \frac{n}{n-i}  \binom{n-i}{i} (xy)^i ((x+y)^2)^{\left\lfloor \frac{n}{2} \right\rfloor-i}.
	\end{equation}
	
	From \eqref{comp1} and \eqref{00}, we get the proof of the following explicit formula:
	
	\begin{theorem}
		\label{comp2}
		For any natural number $n$, the following formula is true:
		\begin{equation}
			\label{comp3}
			\Psi(a,b,n) =\sum_{i=0}^{\left\lfloor \frac{n}{2} \right\rfloor}\frac{n}{n-i} \binom{n-i}{i} (-a)^i (2a-b)^{\left\lfloor \frac{n}{2} \right\rfloor - i}.
		\end{equation}
	\end{theorem}

	\subsection{GENERALIZATION OF THE EIGHT LEVEL THEOREM}	
	\begin{theorem}{(Generalization of the Eight Level Theorem)}\\
		\label{exp1}
		For any natural number $n$, and any real numbers $a,b, \alpha, \beta$, $ \beta a - \alpha b \neq 0 $, there exist unique polynomials in $a,b, \alpha, \beta$ with integer coefficients, that we call
		$ \Psi\left( \begin{array}{cc|r} a & b & n \\ \alpha & \beta & r \end{array} \right)$, that depend only on $a,b, \alpha, \beta, n,$ and $r$, and satisfy the following polynomial expansion:
		\begin{equation}
			\label{000} 
			(\beta a - \alpha b)^{\lfloor{\frac{n}{2}}\rfloor} \frac{x^n+y^n}{(x+y)^{\delta(n)}}  = \sum_{r=0}^{\lfloor{\frac{n}{2}}\rfloor}
			\Psi\left(\begin{array}{cc|r} a & b & n \\ \alpha & \beta & r \end{array}\right)
			(\alpha x^2 + \beta xy + \alpha y^2)^{\lfloor{\frac{n}{2}}\rfloor -r} (ax^2+bxy+ay^{2})^{r}. 
		\end{equation}
		Moreover, we have the following:
		\begin{equation}
			\label{000-a} 
			\Psi\left(\begin{array}{cc|r} a & b & n \\ \alpha & \beta & 0 \end{array} \right) = \Psi(a,b,n),
		\end{equation}
		\begin{equation}
			\label{000-b} 
		\qquad \quad 	\Psi\left(\begin{array}{cc|c} a & b & n \\ \alpha & \beta & \lfloor{\frac{n}{2}}\rfloor \end{array} \right) = (-1)^{\lfloor{\frac{n}{2}}\rfloor} \: \Psi(\alpha,\beta,n). 
		\end{equation}
	\end{theorem}

\begin{proof}
From \eqref{00}, and noting that
\begin{equation}
	\label{bex1}
	\begin{aligned}
		(\beta a-\alpha b)\: (x+y)^2 & = (2a-b)(\alpha x^2 + \beta xy + \alpha y^2) + (\beta- 2\alpha)( a x^2 + b xy + a y^2   ),\\	
		(\beta a-\alpha b) xy & = a (\alpha x^2 + \beta xy + \alpha y^2) + (-\alpha)( a x^2 + b xy + a y^2   ),
	\end{aligned}
\end{equation}
we get the proof of the existence of the  polynomials $\Psi\left(\begin{array}{cc|c} a & b & n \\ \alpha & \beta & r \end{array} \right)$ with integer coefficients. 
The uniqueness of the coefficients come from the fact that $ \alpha x^2 + \beta xy + \alpha y^2 $ and $   a x^2 + b xy + a y^2$ are algebraically independent for $ \beta a - \alpha b \neq 0 $. Put $x=x_0=-b + \sqrt{b^2 - 4a^2}, \:y= y_0=2a$, then $ax_0^2 + bx_0 y_0 + a y_0^2 = 0$. It follows that 
\begin{equation}
	\begin{aligned}
		\label{comp12}	 	
		\Psi\left(\begin{array}{cc|r} a & b & n \\ \alpha & \beta & 0 \end{array} \right) &= \frac{(2a-b)^{\lfloor{\frac{n}{2}}\rfloor}} {2^n} \left\{  \left( 1 + \sqrt{ \frac{b+2a}{b-2a}} \right)^n + \left( 1 - \sqrt{ \frac{b+2a}{b-2a}} \right)^n \right\}\\
		&= \Psi(a,b,n).
	\end{aligned}	 	
\end{equation} 	
From \eqref{comp12} and \eqref{comp1} we get \eqref{000-a}.
Now put $x=x_1=-\beta + \sqrt{\beta^2 - 4 \alpha^2}, \:y= y_1=2 \alpha$, then $\alpha x_1^2 + \beta x_1 y_1 + \alpha y_1^2 = 0$. Hence  
\begin{equation}
	\begin{aligned}
		\label{comp123}	 	
		\Psi\left(\begin{array}{cc|c} a & b & n \\ \alpha & \beta & \lfloor{\frac{n}{2}}\rfloor \end{array} \right) &= (-1)^{\lfloor{\frac{n}{2}}\rfloor} \: \frac{(2 \alpha -\beta)^{\lfloor{\frac{n}{2}}\rfloor}} {2^n} \left\{  \left( 1 + \sqrt{ \frac{\beta+2\alpha}{\beta-2\alpha}} \right)^n + \left( 1 - \sqrt{ \frac{\beta+2\alpha}{\beta-2\alpha}} \right)^n \right\}\\
		&= (-1)^{\lfloor{\frac{n}{2}}\rfloor} \: \Psi(\alpha,\beta,n).
	\end{aligned}	 	
\end{equation} 	
From \eqref{comp123} and \eqref{comp1} we get \eqref{000-b}. This completes the proof.	\end{proof}

\subsection{\textbf{THE FIRST FUNDAMENTAL THEOREM OF $\Psi$-SEQUENCE}}

\begin{theorem}{(The first fundamental theorem of $\Psi$-sequence)}\\
	\label{Aexp1} For $\beta a - \alpha b \neq 0$,
	the polynomials $\Psi_r(n):=\Psi\left(\begin{array}{cc|c}
		a & b & n \\ \alpha & \beta & r \end{array}\right)$ satisfy
	\begin{align}
		\begin{aligned}
			\label{diff1}
			\big(\alpha \frac{{\partial} }{\partial a} +  \beta \frac{{\partial} }{\partial b}\big)\Psi_r(n) &= - (r+1)\Psi_{r+1}(n), \\
			\big(a \frac{{\partial} }{\partial \alpha} +  b \frac{{\partial} }{\partial \beta} \big) \Psi_r(n) &= - \big(\lfloor{\frac{n}{2}}\rfloor - r + 1 \big)\Psi_{r-1}(n).
		\end{aligned}
	\end{align}
\end{theorem}

\begin{proof}
	We differentiate \eqref{00} with respect to the particular differential operator 
	\[ 		\big(\alpha \frac{{\partial} }{\partial a} +  \beta \frac{{\partial} }{\partial b} \big).    \]
	Noting that
	\begin{equation}
		(\alpha \frac{{\partial} }{\partial a} +  \beta \frac{{\partial} }{\partial b}) (\beta a - \alpha b)^{\lfloor{\frac{n}{2}}\rfloor} \frac{x^n+y^n}{(x+y)^{\delta(n)}} =0,
	\end{equation}
	we get
	\begin{equation}
		\label{p1}
		\begin{aligned}
			0=\Big(\alpha \frac{{\partial} }{\partial a} & +  \beta \frac{{\partial} }{\partial b} \Big) \sum_{r=0}^{\lfloor{\frac{n}{2}}\rfloor}
			\Psi_r(n) (\alpha x^2 + \beta xy + \alpha y^2)^{\lfloor{\frac{n}{2}}\rfloor -r} (ax^2+bxy+ay^{2})^{r}  \\
			&=  \sum_{r=0}^{\lfloor{\frac{n}{2}}\rfloor}
			(\alpha x^2 + \beta xy + \alpha y^2)^{\lfloor{\frac{n}{2}}\rfloor -r} (ax^2+bxy+ay^{2})^{r}\Big(\alpha \frac{{\partial} }{\partial a} +  \beta \frac{{\partial} }{\partial b} \Big)\Psi_r(n)\\
			&+   \sum_{r=0}^{\lfloor{\frac{n}{2}}\rfloor}
			\Psi_r(n)
			(\alpha x^2 + \beta xy + \alpha y^2)^{\lfloor{\frac{n}{2}}\rfloor -r}\Big(\alpha \frac{{\partial} }{\partial a} +  \beta \frac{{\partial} }{\partial b} \Big) (ax^2+bxy+ay^{2})^{r} \\
			&+  \sum_{r=0}^{\lfloor{\frac{n}{2}}\rfloor}
			\Psi_r(n)
			(ax^2+bxy+ay^{2})^{r}\Big(\alpha \frac{{\partial} }{\partial a} +  \beta \frac{{\partial} }{\partial b} \Big) (\alpha x^2 + \beta xy + \alpha y^2)^{\lfloor{\frac{n}{2}}\rfloor -r}.
		\end{aligned}
	\end{equation}
	Consequently, from \eqref{p1}, we obtain the following desirable polynomial expansion
	\begin{equation}
		\label{p2}
		0 = \sum_{r=0}^{\lfloor{\frac{n}{2}}\rfloor} \Big(	(\alpha \frac{{\partial} }{\partial a} +  \beta \frac{{\partial} }{\partial b})\Psi_r(n) + (r+1)\Psi_{r+1}(n) \Big)  (\alpha x^2 + \beta xy + \alpha y^2)^{\lfloor{\frac{n}{2}}\rfloor -r} (ax^2+bxy+ay^{2})^{r}.
	\end{equation}
	As $ \beta a - \alpha b \neq 0 $, the polynomials $(\alpha x^2 + \beta xy + \alpha y^2)$ and $(a x^2 + b xy + a y^2)$ are algebraically independent which means that all of the coefficients of ~\eqref{p2} must vanish. This means that 
	\begin{align}
		\begin{aligned}
			\big(\alpha \frac{{\partial} }{\partial a} +  \beta \frac{{\partial} }{\partial b}\big)\Psi_r(n) + (r+1)\Psi_{r+1}(n)  = 0 \quad   \text{for all} \:\: r.
		\end{aligned}
	\end{align}
	Similarly, we differentiate \eqref{00} with respect to the particular differential operator 
	\[ 		\big(a \frac{{\partial} }{\partial \alpha} +  b \frac{{\partial} }{\partial \beta} \big).    \]
	We get
	\begin{align}
		\begin{aligned}
			\big(a \frac{{\partial} }{\partial \alpha} +  b \frac{{\partial} }{\partial \beta} \big) \Psi_r(n) + \big(\lfloor{\frac{n}{2}}\rfloor - r + 1 \big)\Psi_{r-1}(n) = 0  \quad   \text{for all} \:\: r. 
		\end{aligned}
	\end{align}
	This completes the proof.
\end{proof}
This result immediately gives the following desirable theorem

\begin{theorem}
	\label{A1} For $\beta a - \alpha b \neq 0$,
	the polynomials $\Psi\left(\begin{array}{cc|c}
		a & b & n \\ \alpha & \beta & r \end{array}\right)$ satisfy
	\begin{align}
		\begin{aligned}
			\label{diff3}
			\Psi\left(\begin{array}{cc|c}
				a & b & n \\ \alpha & \beta & r \end{array} \right) &=  \frac{(-1)^r}{r!}
			\Big(\alpha \frac{{\partial} }{\partial a} + \beta \frac{{\partial}}{\partial b}\Big)^{r} \Psi(a,b,n),  \quad \qquad
		\end{aligned}
	\end{align}
	and
	\begin{align}
		\begin{aligned}
			\label{diff5}
			\Psi\left( \begin{array}{cc|c}
				a & b & n \\ \alpha & \beta & r \end{array} \right) =  \frac{(-1)^r}{(\lfloor{\frac{n}{2}}\rfloor -r)!} \Big(a \frac{{\partial} }{\partial \alpha} + b \frac{{\partial}}{\partial \beta}\Big)^{\lfloor{\frac{n}{2}}\rfloor-r} \Psi(\alpha,\beta,n).  
		\end{aligned}
	\end{align}
\end{theorem}

\section*{\textbf{EXPLORING TWO ILLUSTRATIVE EXAMPLES IN DETAIL}}

\subsection*{\textbf{ILLUSTRATIVE EXAMPLE - 1}}
To illustrate a simple example for Theorem \eqref{exp1}, Theorem \eqref{A1},  let's consider $n=4$. Simple computations yield:
\begin{equation}
	\begin{aligned}
		\Psi\left(\begin{array}{cc|r} a & b & n \\ \alpha & \beta & 0 \end{array}\right)	&= \frac{(-1)^0}{0!}
		\Big(\alpha \frac{{\partial} }{\partial a} + \beta \frac{{\partial}}{\partial b}\Big)^{0} \Psi(a,b,4) = \Psi(a,b,4) = -2 a^2 + b^2 , \\
		\Psi\left(\begin{array}{cc|r} a & b & n \\ \alpha & \beta & 1 \end{array}\right) &=\frac{(-1)^1}{1!}
		\Big(\alpha \frac{{\partial} }{\partial a} + \beta \frac{{\partial}}{\partial b}\Big)^{1} \Psi(a,b,4) = 4a \alpha - 2b \beta , \\
		\Psi\left(\begin{array}{cc|r} a & b & n \\ \alpha & \beta & 2 \end{array}\right) &=\frac{(-1)^2}{2!}
		\Big(\alpha \frac{{\partial} }{\partial a} + \beta \frac{{\partial}}{\partial b}\Big)^{2} \Psi(a,b,4) = -2 \alpha^2 + \beta^2.
	\end{aligned}
\end{equation}

Hence, we derive the following polynomial identity, which represents a specific instance of the primary results outlined in this paper:
\begin{align}
	\begin{aligned}
		\label{example}
		(\beta a - \alpha b)^2 (x^4 + y^4) &= (-2 a^2 + b^2)(\alpha x^2 +\beta xy + \alpha y^2)^2 \\
		&+ (4a \alpha - 2b \beta)(\alpha x^2 +\beta xy + \alpha y^2)(a x^2 +b xy + a y^2) \\
		&+ (-2 \alpha^2 + \beta^2)(a x^2 + b xy + a y^2)^2.
	\end{aligned}
\end{align}

This example not only provides an intuitive understanding of Equation \eqref{Main theorem 1} but also reveals a fascinating generalization of a renowned identity in number theory. For the case when $n=4$ and $\left( \begin{array}{cc} a & b \\ \alpha & \beta \end{array} \right) = \left( \begin{array}{cc} 1 & 1 \\ 1 & 2 \end{array} \right)$, the middle coefficient $4a \alpha - 2b \beta$ in \eqref{Main theorem 1} becomes zero, leading us to a special instance of a well-established identity often used in the exploration of equal sums of like powers and the derivation of new formulas for Fibonacci numbers:

\begin{align}
	\begin{aligned}
		x^4 + y^4 + (x+y)^4 = 2 (x^2 + xy + y^2)^2.
	\end{aligned}
\end{align}

In Volume 2, on page 650, \cite{Dickson} attributes this special case to C. B. Haldeman (1905), although Proth (1878) briefly mentioned it earlier (see page 657 of \cite{Dickson}).

\subsection*{\textbf{ILLUSTRATIVE EXAMPLE - 2}}

Take $n=6$. We know that 
\[ \Psi(a,b,6)= 3a^2b-b^3. \]  
Therefore, 
\begin{align}
	\begin{aligned}
	(\frac{{\partial} }{\partial a} + 2 \frac{{\partial}}{\partial b})(3a^2b-b^3) &= 6a^2+6ab-6b^2 , \\
	(\frac{{\partial} }{\partial a} + 2 \frac{{\partial}}{\partial b})(6a^2+6ab-6b^2) &= 24a-18b, \\
	(\frac{{\partial} }{\partial a} + 2 \frac{{\partial}}{\partial b})(24a-18b) &= -12. \\
	\end{aligned}
\end{align}
Therefore, we derive the following expansion from the Generalized Eight Levels theorem and the first fundamental theorem of the $\Psi$-sequence:
\begin{align}
	\begin{aligned}
		(2a-b)^3(x^6+y^6) &= (3a^2b-b^3) (x+y)^6 \\ &+ (-6a^2-6ab+6b^2) (x+y)^4 (ax^2 +bxy+ay^2) \\ &+  (12a-9b) (x+y)^2(ax^2 +bxy+ay^2)^2 + (2)(ax^2 +bxy+ay^2)^3.
	\end{aligned}
\end{align}

\section*{\textbf{EXPLORING TWO EXPLICIT FORMULAS}}

\subsection*{\textbf{EXPLICIT FORMULA - 1}}

Upon reviewing the expansion in \eqref{000} of the Generalized Eight Levels Theorem \eqref{exp1} and equation \eqref{00}, it becomes apparent that the subsequent formula is valid for all integer values of $n$:

\begin{equation}
	\label{formula 1}
	\Psi\left( \begin{array}{cc|r} 0 & 1 & n \\ 1 & 2 & r \end{array} \right) = (-1)^{\lfloor{\frac{n}{2}}\rfloor-r} \frac{n}{n-r} \binom{n-r}{r}.
\end{equation}

\subsection*{\textbf{EXPLICIT FORMULA - 2}}

It is worth noting that the following result can be derived using the binomial theorem. We can express it as:

\begin{equation}
	\label{special sum}
	\begin{aligned}
		4^{\lfloor{\frac{n}{2}}\rfloor}	\frac{x^n+y^n}{(x+y)^{\delta(n)}} = 
		\sum_{r=0}^{\lfloor{\frac{n}{2}}\rfloor}  2^{\delta(n-1)} \binom{n}{2r} (x^2 +2xy+y^2 )^{\lfloor{\frac{n}{2}}\rfloor -r} (x^2-2xy+y^{2})^{r}.
	\end{aligned}
\end{equation}

Hence, we can conclude:

\begin{equation}
	\label{formula 2}
	\Psi\left( \begin{array}{cc|r} 1 & -2 & n \\ 1 & 2 & r \end{array} \right) = 2^{\delta(n-1)}  \binom{n}{2r}.
\end{equation}

\subsection{\textbf{THE SECOND FUNDAMENTAL THEOREM OF $\Psi$-SEQUENCE}}
Now, put $ r=\lfloor{\frac{n}{2}}\rfloor  $ in equation \eqref{diff3} of Theorem \eqref{A1}, together with Theorem \eqref{exp1}, we get the following immediate consequence

\begin{theorem}{(The second fundamental theorem of $\Psi$-sequence)}\\
	\label{Aexp2} For any numbers $a,b, \alpha, \beta  $, $\beta a - \alpha b \neq 0$, and any natural number $n$, we have
	\begin{align}
		\begin{aligned}
			\frac{1}{(\lfloor{\frac{n}{2}}\rfloor)!} \Big(\alpha \frac{{\partial} }{\partial a} + \beta \frac{{\partial}}{\partial b}\Big)^{\lfloor{\frac{n}{2}}\rfloor} \Psi(a,b,n) = \Psi(\alpha,\beta,n).  
		\end{aligned}
	\end{align}
\end{theorem}		
Also, it is useful to deduce the following relations. Replace each $\alpha$ and $\beta$ by $\lambda \alpha$ and $\lambda \beta$ respectively in Theorem \eqref{exp1}, we get the following polynomial identity for any $\lambda$
\begin{equation*}
	\label{}
	\begin{aligned}
		&(\lambda\beta a - \lambda\alpha b)^{\lfloor{\frac{n}{2}}\rfloor} \frac{x^n+y^n}{(x+y)^{\delta(n)}} = \\
		&\sum_{r=0}^{\lfloor{\frac{n}{2}}\rfloor}
		\Psi\left( \begin{array}{cc|r} a & b & n \\ \lambda\alpha & \lambda\beta & r \end{array} \right)
		(\lambda\alpha x^2 + \lambda\beta xy + \lambda \alpha y^2)^{\lfloor{\frac{n}{2}}\rfloor -r} (ax^2+bxy+ay^{2})^{r}.
	\end{aligned}
\end{equation*}
Then
\begin{equation}
	\label{V1}
	\begin{aligned}
		&(\beta a - \alpha b)^{\lfloor{\frac{n}{2}}\rfloor} \frac{x^n+y^n}{(x+y)^{\delta(n)}} = \\ &\sum_{r=0}^{\lfloor{\frac{n}{2}}\rfloor} \lambda^{ -r}
		\Psi\left( \begin{array}{cc|r} a & b & n \\ \lambda\alpha & \lambda\beta & r \end{array} \right)
		(\alpha x^2 + \beta xy +  \alpha y^2)^{\lfloor{\frac{n}{2}}\rfloor -r} (ax^2+bxy+ay^{2})^{r}.
	\end{aligned}
\end{equation}
Comparing \eqref{V1} with \eqref{exp1}, we obtain
\[     \lambda^{ -r}
\Psi\left( \begin{array}{cc|r} a & b & n \\ \lambda\alpha & \lambda\beta & r \end{array} \right) =     \Psi\left( \begin{array}{cc|r} a & b & n \\ \alpha & \beta & r \end{array} \right).       \]
Similarly, we can prove the following useful relations.

\begin{theorem}
	\label{W11}
	For any numbers $a, b, \alpha, \beta, \beta a - \alpha b \neq 0, \lambda, r,  n $,  we get
	\begin{equation}
		\label{W22}
		\begin{aligned}
			\Psi\left( \begin{array}{cc|c}
				a & b & n \\ \lambda \alpha & \lambda \beta & r \end{array} \right)
			&=\lambda^{r} \Psi\left( \begin{array}{cc|c}
				a & b & n \\ \alpha & \beta & r \end{array} \right),  \\  
			\Psi\left( \begin{array}{cc|c}
				\lambda a &\lambda b & n \\ \alpha & \beta & r \end{array} \right)
			&=\lambda^{\lfloor{\frac{n}{2}}\rfloor - r}
			\Psi\left(\begin{array}{cc|c}
				a & b & n \\ \alpha & \beta & r \end{array} \right), \\
			\Psi\left( \begin{array}{cc|c}
				a & b & n \\ \alpha & \beta & r \end{array} \right) &=
			(-1)^{\lfloor{\frac{n}{2}} \rfloor}
			\Psi\left( \begin{array}{cc|c}
				\alpha & \beta  & n \\  a & b &\lfloor{\frac{n}{2}} \rfloor - r \end{array} \right),
			\\
			\text{and}\quad \quad \quad \quad  \quad \quad \quad \quad \quad \quad \quad \quad \\
			\lambda^{\lfloor{\frac{n}{2}}\rfloor}\Psi(a,b,n) &= \Psi(\lambda a,\lambda  b,n).  		
		\end{aligned}
	\end{equation}
\end{theorem}

\section{\textbf{REPRESENTATIONS FOR THE $\Psi$-SEQUENCE}}
We now ready to prove the following theorem.
\begin{theorem}
	\label{Ready00}
	For any $a,b,\alpha,\beta, \theta, n$, $\beta a - \alpha b \neq 0,$ the following identities are true
	\begin{equation}
		\label{Ready1}
		\sum_{r=0}^{\lfloor{\frac{n}{2}}\rfloor}  \Psi\left( \begin{array}{cc|r} a & b & n \\ \alpha & \beta & r \end{array} \right) \theta^r = \Psi(a-\alpha \theta , b - \beta \theta, n).
	\end{equation}
	
\end{theorem}
\begin{proof}
	Define $q_1:=\alpha x^2+\beta xy +\alpha y^2$ and $q_2:=a x^2+b xy +a y^2$ and
	\[ \Lambda_\theta :=\theta q_1 - q_2 = (\alpha \theta - a)x^2 + (\beta \theta - b)xy + (\alpha \theta - a)y^2 .\] 	
	From Theorem \eqref{exp1}, we know that
	\[
	(\beta a - \alpha b)^{\lfloor{\frac{n}{2}}\rfloor} \frac{x^n+y^n}{(x+y)^{\delta(n)}} = \sum_{r=0}^{\lfloor{\frac{n}{2}}\rfloor}  \Psi\left( \begin{array}{cc|r} a & b & n \\ \alpha & \beta & r \end{array} \right) (q_1)^{\lfloor{\frac{n}{2}}\rfloor -r} (q_2)^{r}. \]
	As $q_2 \equiv \theta q_1  \pmod{\Lambda_\theta}$, we get
	\begin{equation}
		\label{Ready2}
		(\beta a - \alpha b)^{\lfloor{\frac{n}{2}}\rfloor} \frac{x^n+y^n}{(x+y)^{\delta(n)}} \equiv q_1^{\lfloor{\frac{n}{2}}\rfloor} \sum_{r=0}^{\lfloor{\frac{n}{2}}\rfloor}  \Psi\left( \begin{array}{cc|r} a & b & n \\ \alpha & \beta & r \end{array} \right) \theta^r   \pmod{\Lambda_\theta}.
	\end{equation}
	Replace each of $a,b$ by $\alpha \theta - a, \beta \theta - b$ respectively,  in Theorem \eqref{exp1}, we obtain
	\begin{equation}
		\label{}
		(\beta [\alpha \theta - a] -
		\alpha[\beta \theta - b])^{\lfloor{\frac{n}{2}}\rfloor} \frac{x^n+y^n}{(x+y)^{\delta(n)}} \equiv
		\Psi(\alpha \theta - a, \beta \theta - b, n)  q_1^{\lfloor{\frac{n}{2}}\rfloor} \pmod{\Lambda_\theta}.
	\end{equation}	
	As $\beta [\alpha \theta - a] -
	\alpha[\beta \theta - b] = - (\beta a - \alpha b)$, and noting from Theorem \eqref{W11} that
	\[ (-1)^{\lfloor{\frac{n}{2}}\rfloor} \Psi(\alpha \theta - a, \beta \theta - b, n) =   \Psi(a-\alpha \theta , b - \beta \theta, n), \]
	we immediately get the following congruence
	\begin{equation}
		\label{Ready3}
		(  \beta a - \alpha b  )^{\lfloor{\frac{n}{2}}\rfloor} \frac{x^n+y^n}{(x+y)^{\delta(n)}} \equiv
		\Psi(a-\alpha \theta , b - \beta \theta, n) q_1^{\lfloor{\frac{n}{2}}\rfloor} \pmod{\Lambda_\theta}.
	\end{equation}
	Now, subtracting \eqref{Ready2} and\eqref{Ready3}, we obtain
	\begin{equation}
		\label{Ready4}
		0 \equiv \Big(
		\sum_{r=0}^{\lfloor{\frac{n}{2}}\rfloor}  \Psi\left( \begin{array}{cc|r} a & b & n \\ \alpha & \beta & r \end{array} \right) \theta^r -
		\Psi(a-\alpha \theta , b - \beta \theta, n) \Big) q_1^{\lfloor{\frac{n}{2}}\rfloor} \pmod{\Lambda_\theta}.
	\end{equation}
	As the congruence \eqref{Ready4} is true for any $x,y$, and as $(\beta \theta - b) \alpha -  (\alpha \theta - a) \beta = \beta a - \alpha b \neq 0$, then the binary quadratic forms $\Lambda_\theta$ and $ q_1 $ are algebraically independent. This immediately leads to
	\begin{equation}
		\label{Ready5}
		0= \sum_{r=0}^{\lfloor{\frac{n}{2}}\rfloor}  \Psi\left( \begin{array}{cc|r} a & b & n \\ \alpha & \beta & r \end{array} \right) \theta^r -
		\Psi(a-\alpha \theta , b - \beta \theta, n).
	\end{equation}
	Hence we obtained the proof of \eqref{Ready1}. This completes the proof of Theorem \eqref{Ready00}.  \end{proof}

\section{\textbf{SPECIALIZATION AND LIFTING}}
The following desirable generalization is important

\begin{theorem}
	\label{Ready6}
	For any $a,b,\alpha,\beta, \eta, \xi, n$, $\beta a - \alpha b \neq 0,$ the following identities are true
	\begin{equation}
		\label{Ready66}
		\begin{aligned}
			\sum_{r=0}^{\lfloor{\frac{n}{2}}\rfloor}  \Psi\left( \begin{array}{cc|r} a & b & n \\ \alpha & \beta & r \end{array} \right)  \xi^{\lfloor{\frac{n}{2}}\rfloor  - r}  \eta^r =\Psi(a \xi-\alpha \eta, b \xi - \beta \eta, n).
		\end{aligned}
	\end{equation}
\end{theorem}
\begin{proof}
	Without loss of generality, let $\xi \neq 0$. We obtain the proof by replacing each $\theta$ in equation \eqref{Ready1} of Theorem \eqref{Ready00} by $\frac{\eta}{\xi},$ and multiplying each side by $\xi^{\lfloor{\frac{n}{2}}\rfloor},$ and  noting from Theorem \eqref{W11} that 
	\[\xi^{\lfloor{\frac{n}{2}}\rfloor}\Psi(a-\alpha \frac{\eta}{\xi} , b - \beta \frac{\eta}{\xi}, n) = \Psi(a \xi-\alpha \eta, b \xi - \beta \eta, n). \] 
\end{proof}
Replacing $\theta$ by $\pm 1$ in \eqref{Ready1}, we obtain the following desirable special cases
\begin{theorem}
	\label{sum}
	For any $a,b,\alpha,\beta, n$, $\beta a - \alpha b \neq 0$ the following identities are true
	\begin{align}
		\begin{aligned}
			\label{theta equal 1}
			\centering	 	
			\sum_{r=0}^{\lfloor{\frac{n}{2}}\rfloor}  \Psi\left( \begin{array}{cc|r} a & b & n \\ \alpha & \beta & r \end{array} \right) &= \Psi(a-\alpha , b - \beta , n), \\ 	 	
			\sum_{r=0}^{\lfloor{\frac{n}{2}}\rfloor}  \Psi\left( \begin{array}{cc|r} a & b & n \\ \alpha & \beta & r \end{array} \right) (-1)^{r} &= \Psi(a + \alpha , b + \beta , n). 
		\end{aligned}
	\end{align}
\end{theorem}

\subsection{MORE GENERALIZATIONS}
Now we can generalize Theorem \eqref{Ready00} by applying the following specific differential map
\[ \Big( - \frac{{\partial}}{\partial \theta} \Big)^{k}   \]
on equation \eqref{Ready1}, and  noting that 

\begin{equation*}
	\begin{aligned}
		\Big( - \frac{{\partial}}{\partial \theta} \Big)^{k} &
		\Psi(a-\alpha \theta , b - \beta \theta, n) 
		= (-1)^{k} (k!) \Psi\left( \begin{array}{cc|r} a-\alpha \theta & b - \beta \theta & n \\ \alpha & \beta & k\end{array} \right).	
	\end{aligned}
\end{equation*}
Hence we immediately obtain the following desirable generalization
\begin{theorem}
	\label{Ready7}
	For any $n,k, a,b,\alpha,\beta, \theta, n$, $\beta a - \alpha b \neq 0,$ the following identity is true
	\begin{equation}
		\label{Ready8}
		\sum_{r=k}^{\lfloor{\frac{n}{2}}\rfloor} \binom{r}{k}  \Psi\left( \begin{array}{cc|r} a & b & n \\ \alpha & \beta & r \end{array} \right) \theta^{r-k} = \Psi\left( \begin{array}{cc|r} a-\alpha \theta  & b-\beta \theta  & n \\ \alpha & \beta & k \end{array} \right).
	\end{equation}
	
\end{theorem}
Again, without loss of generality, let $\xi \neq 0$. By replacing each $\theta$ in \eqref{Ready8} by $\frac{\eta}{\xi}$, and multiplying each side by $\xi^{\lfloor{\frac{n}{2}}\rfloor -k}$, and noting the properties of $\Psi$ of Theorem \eqref{W11} that

\begin{equation*} \xi^{\lfloor{\frac{n}{2}}\rfloor -k} \Psi\left( \begin{array}{cc|r} a-\alpha \frac{\eta}{\xi}  & b-\beta \frac{\eta}{\xi}  & n \\ \alpha & \beta & k \end{array}     \right)  =  \Psi\left( \begin{array}{cc|r} a\xi-\alpha \eta  & b\xi-\beta \eta  & n \\ \alpha & \beta & k \end{array} \right),        		\end{equation*}

we obtain the following generalization for Theorem \eqref{Ready7}
\begin{theorem}
	\label{Ready9}
	For any $n,k,a,b,\alpha,\beta, \theta, n$, $\beta a - \alpha b \neq 0,$ the following identity is true
	\begin{align}
		\begin{aligned}
			\label{Ready10}
			\sum_{r=k}^{\lfloor{\frac{n}{2}}\rfloor} \binom{r}{k}  \Psi\left( \begin{array}{cc|r} a & b & n \\ \alpha & \beta & r \end{array} \right) \xi^{\lfloor{\frac{n}{2}}\rfloor  - r}  \eta^{r-k} = \Psi\left( \begin{array}{cc|r} a\xi-\alpha \eta  & b\xi-\beta \eta  & n \\ \alpha & \beta & k \end{array} \right).
		\end{aligned}
	\end{align}
\end{theorem}
\section{\textbf{THE $\Psi$-REPRESENTATION FOR THE SUM OF LIKE POWERS}}
Now, put $\xi = \alpha x^2 +\beta xy + \alpha y^2$ and $ \eta = a x^2 +bxy +a y^2$ in equation \eqref{Ready66} of Theorem \eqref{Ready6}, we obtain
\begin{equation}
	\begin{aligned}
		\label{WW1}
		\sum_{r=0}^{\lfloor{\frac{n}{2}}\rfloor}  \Psi\left( \begin{array}{cc|r} a & b & n \\ \alpha & \beta & r \end{array} \right) & (\alpha x^2 +\beta xy + \alpha y^2)^{\lfloor{\frac{n}{2}}\rfloor  - r}  (  a x^2 +bxy +a y^2 )^r \\
		& = (  \beta a - \alpha b  )^{\lfloor{\frac{n}{2}}\rfloor} \Psi(xy,-x^2-y^2, n).
	\end{aligned}
\end{equation}
Now, from \eqref{exp1}, we get
\begin{equation}
	\begin{aligned}
		\label{WW2}
		\sum_{r=0}^{\lfloor{\frac{n}{2}}\rfloor}  \Psi\left( \begin{array}{cc|r} a & b & n \\ \alpha & \beta & r \end{array} \right) & (\alpha x^2 +\beta xy + \alpha y^2)^{\lfloor{\frac{n}{2}}\rfloor  - r}  (  a x^2 +bxy +a y^2 )^r \\
		&= (  \beta a - \alpha b  )^{\lfloor{\frac{n}{2}}\rfloor}\frac{x^n+y^n}{(x+y)^{\delta(n)}}.
	\end{aligned}
\end{equation}
From \eqref{WW1}, \eqref{WW2}, we get the following desirable $\Psi-$representation for the sums of powers 	
\begin{theorem}{(The $\Psi-$representation for sums of powers)}
	\label{WW3}
	For any natural number $n$, the $\Psi-$polynomial satisfy the following identity
	\begin{equation}
		\label{WW4}
		\begin{aligned}
			\Psi(xy,-x^2-y^2,n) &= \frac{x^n+y^n}{(x+y)^{\delta(n)}}.   \\
		\end{aligned}
	\end{equation}
\end{theorem}

\section{\textbf{NEW FORMULAS FOR SUM OF LIKE POWERS}}

From Theorem \eqref{WW3}, we can readily establish that:
\begin{align}
	\begin{aligned}
		\Psi(zt,-(z^2+t^2),n) &= \frac{z^n+t^n}{(z+t)^{\delta(n)}}, \\
		\Psi(uv,-(u^2+v^2),n) &= \frac{u^n+v^n}{(u+v)^{\delta(n)}}.  
	\end{aligned}
\end{align}

It becomes apparent that, according to the Generalized Eight Levels theorem \eqref{exp1}, the first and last coefficients on the right-hand side of \eqref{000} are $\Psi(a,b,n)$ and $(-1)^{\lfloor{\frac{n}{2}}\rfloor} \Psi(\alpha,\beta,n)$, respectively. This observation, along with the intriguing natural symmetry, motivates us to investigate a specific special case of the Generalized Eight Levels theorem. To do so, we consider the substitution of a new set of parameters $z, t, u, v$:
\[
a = zt,  b =-(z^2+t^2), \alpha = uv, \beta = - (u^2+v^2). 
\] 
For the sake of simplification, let's introduce the following notation:

\begin{equation}
	\label{n=notation}
	\boldmath
	\begin{aligned}
		\centering
			[x,y \vert u,v ] := (xu - yv) (xv - yu) = (x^2+y^2)uv-xy(u^2+v^2) 
	\end{aligned}
\end{equation}

Observe that this symbol exhibits the following elementary properties:

\begin{align*}
	[x,y \vert u,v] &= - [u,v \vert x,y ] &  [dx,dy \vert u,v ] &= d^2 [x,y \vert u,v ] & [x,y \vert du,dv ] &= d^2 [x,y \vert u,v ]    \\
	[x,y \vert u,v ] &= [y,x \vert u,v ] &
	[x,y \vert u,v ] &= [x,y \vert v,u ]  & [x,y \vert u,v ] &= [y,x \vert v,u]    \\
	[0,1 \vert u,v ] &= uv & [1,1 \vert u,v ] &= - (u-v)^2 & 
	[1,-1 \vert u,v ] &= (u+v)^2 \\
	[1,-i \vert u,v ] &= i [u^2+v^2] &
	[1,i \vert u,v ] &= -i (u^2+v^2) & [i,i \vert u,v ] &= (u-v)^2     
\end{align*}
Upon substituting \[ \Psi\left( \begin{array}{cc|r} a & b & n \\ \alpha & \beta & r \end{array} \right) = \Psi\left( \begin{array}{cc|r} zt & z^2+t^2 & n \\ uv & u^2+v^2 & r \end{array} \right) \]
into the Generalized Eight Levels theorem, we promptly arrive at the following noteworthy special case:

\begin{theorem}{(Special Case of the Generalized Eight Levels Theorem)}
\label{special6}

	The following polynomial expansion is true
	\begin{equation}
		\label{special-5}
		\begin{aligned}
			&\left[z,t \vert u,v\right]^{\lfloor{\frac{n}{2}}\rfloor} \frac{x^n+y^n}{(x+y)^{\delta(n)}} - 
			\left[x,y \vert u,v \right]^{\lfloor{\frac{n}{2}}\rfloor} 
			\frac{z^n+t^n}{(z+t)^{\delta(n)}}  -  
			[z,t \vert x,y]^{\lfloor{\frac{n}{2}}\rfloor}
			\frac{u^n+v^n}{(u+v)^{\delta(n)}}  \\
			&= \sum_{r=1}^{\lfloor{\frac{n}{2}}\rfloor-1} 
			\frac{1}{r!} [x,y \vert u,v ]^{\lfloor{\frac{n}{2}}\rfloor -r} [z,t \vert x,y ]^{r} 
			\Big(uv \frac{{\partial} }{\partial zt} + (u^2+v^2) \frac{{\partial}}{\partial z^2+t^2}\Big)^{r} \frac{z^n+t^n}{(z+t)^{\delta(n)}}.
		\end{aligned}
	\end{equation}
	
\end{theorem}
Let's consider some special cases to gain a better understanding of the result.

\section{\textbf{NEW IDENTITIES FOR SUM OF LIKE POWERS}}

\textbf{\large For n=2}

\begin{equation}
	\label{n=2}
	\boldmath
	\begin{aligned}
		\centering
		\left[z,t \vert u,v\right] (x^2+y^2) + \left[u,v \vert x,y \right] (z^2+t^2) +  
		[x,y \vert z,t] (u^2 +v^2) = 0 
	\end{aligned}
\end{equation}

\textbf{\large For n=3}
\begin{equation}
	\label{n=3}
	\boldmath
	\begin{aligned}
		\centering
		\left[z,t \vert u,v\right] \frac{x^3+y^3}{x+y} + \left[u,v \vert x,y \right] \frac{z^3+t^3}{z+t} +  [x,y \vert z,t] \frac{u^3 +v^3}{u+v} = 0 
	\end{aligned}
\end{equation}

\textbf{\large For n=4}

\begin{equation}
	\label{n=4}
	\boldmath
	\begin{aligned}
		\centering
		\left[z,t \vert u,v\right]^2 & (x^4+y^4) - \left[u,v \vert x,y \right]^2 (z^4+t^4) - [x,y \vert z,t]^2 (u^4 +v^4) \\
		& = 2 [x,y \vert u,v] [z,t \vert x,y] ((z^2+t^2) (u^2+v^2)-2 ztuv)
	\end{aligned}
\end{equation}

\textbf{\large For n=5}
\begin{equation}
	\label{fifth power}
	\boldmath
	\begin{aligned}
	&	\left[z,t \vert u,v \right]^{2} \frac{x^5+y^5}{x+y} - 
		[ u,v \vert x,y]^{2} 
		\frac{z^5+t^5}{z+t}  -  
		[x,y \vert z,t]^{2}
		\frac{u^5+v^5}{u+v}  \\
		= [x,y \vert u,v ] &[z,t \vert x,y ] [(u^2+v^2 - uv) (z^2+ t^2) + (z^2+t^2 - zt) (u^2+ v^2) - 2(zt)(uv)]
	\end{aligned}
\end{equation}

\subsection*{REMARKS}
\begin{enumerate}
	\item To compute the expression
	\begin{equation}
		\label{EE}
		\begin{aligned}
			&\left(uv \frac{{\partial} }{\partial zt} + (u^2+v^2) \frac{{\partial}}{\partial z^2+t^2}\right)^{r} \frac{z^n+t^n}{(z+t)^{\delta(n)}},
		\end{aligned}
	\end{equation}
	we need to represent \[ \frac{z^n+t^n}{(z+t)^{\delta(n)}},\]
	as an expansion in terms of the symmetric polynomials \(zt\) and \(z^2+t^2\). The most direct approach to achieve this representation is to apply The Eight Levels Theorem as referenced in \cite{2} or make use of Theorems \eqref{Expansion-3} and \eqref{Expansion-5}:
	
	\begin{equation}
		\label{special-1}
		\begin{aligned}
		\frac{z^n+t^n}{(z+t)^{\delta(n)}} =
		\sum_{k=0}^{\lfloor{\frac{n}{2}}\rfloor}\Psi\left( \begin{array}{cc|r} 1 & 0 & n \\ 0 & 1 & k \end{array} \right) 
		(zt)^{\lfloor{\frac{n}{2}}\rfloor -k} (z^2+t^{2})^{k}.                                              
		\end{aligned}
	\end{equation}

	\item 

From \eqref{n=2} and \eqref{n=3} we get 
\begin{equation}
	\label{n=2-3}
	\boldmath
	\begin{aligned}
		\centering
		\left[z,t \vert u,v\right] (xy) + \left[u,v \vert x,y \right] (zt) +  
		[x,y \vert z,t] (uv) = 0. 
	\end{aligned}
\end{equation}

\item 
As observed, Expansion \eqref{fifth power} represents a specific instance and a unique case of the Expansion \eqref{special-5} from Theorem \eqref{special6}. In 2004, \cite{1}, on page 446, made reference to Expansion \eqref{fifth power} but did not provide any proof. 

\item  Notably, \cite{1} pointed out that for the following values of \(x\), \(y\), \(z\), and \(t\):
\[
x = 5pq, \quad y = 5(p^2 + pq + q^2), \quad z=-5p(p + q), \quad t = -5q(p + q),
\]
we can establish the relationship:
\[
x^5 + y^5 + z^5 + t^5 = d^2,
\]
where \(d\) is given by:
\[
d = 125pq(p + q)(p^2 + pq + q^2).
\]

Based on the structure of \(x\), \(y\), \(z\), \(t\), and \(d\) described above, it is impossible to eliminate the common factor \(5\) in order to find a solution to the equation:
\[
x^5 + y^5 + z^5 + t^5 = d^2
\]
in co-prime polynomials with integer coefficients and without constant common factors. 

\item  \cite{Gawron} succeeded in finding relatively prime solutions and demonstrated that, for \(n\) as a non-zero rational number, the Diophantine equation:
\[
T^2 = n(X_1^5 + X_2^5 + X_3^5 + X_4^5)
\]
has a solution in co-prime polynomials:
\[
X_1, X_2, X_3, X_4, T \in \mathbb{Z}[u, v].
\]

\end{enumerate}

\section{\textbf{FORMULA FOR THE PRODUCTS OF $\Psi(a,b,n) \cdot \Psi(a,b,m) $}}

First, we need to extend the definition of $\Psi(a,b,-l):=\Psi(a,b,l)$ for any integer $l$. Now, for any natural numbers $n,m$, we need a formula for the product $\Psi(a,b,n) \cdot \Psi(a,b,m)$ as a combination of $\Psi(a,b,n+m)$ and $\Psi(a,b,n-m)$. We should observe that the formula \eqref{WW4} is helpful to find that formula. In \eqref{WW4}, put
\begin{equation}
	\begin{aligned}
		x = x_0 = \frac{1}{2}\Big(1 + \sqrt{\frac{b+2a}{b-2a}} \Big),  \\
		y = y_0 = \frac{1}{2}\Big(1 - \sqrt{\frac{b+2a}{b-2a}} \Big).  \\
	\end{aligned}
\end{equation}
Clearly 
	\begin{equation}
	\label{WW-A}
	\begin{aligned}
		x_{0} \; y_{0} = \frac{a}{2a-b} \quad, \quad -x^2_{0} - y^2_{0} = \frac{b}{2a-b} \quad, \quad x_0 +y_0 =1
	\end{aligned}
\end{equation}

 Consequently, from Theorem \eqref{WW3} and Theorem \eqref{W11}, we get 
	\begin{equation}
	\label{WW5}
	\begin{aligned}
		\Psi(a,b,n) = (2a-b)^{\lfloor{\frac{n}{2}}\rfloor} (x_0^n+y_0^n) \quad, \quad  \Psi(a,b,m) =(2a-b)^{\lfloor{\frac{m}{2}}\rfloor} (x_0^m+y_0^m).
	\end{aligned}
\end{equation}
From \eqref{WW5}, and noting the unexpected formula
	\begin{equation}
	\label{WW6}
\lfloor{\frac{n}{2}}\rfloor + \lfloor{\frac{m}{2}}\rfloor - \lfloor{\frac{n+m}{2}}\rfloor + \delta(n) \cdot \delta(m) = 0,
\end{equation}

we immediately get the proof for the following desirable theorem

\begin{theorem}{(The product of $\Psi$-sequences)} \\
	\label{WW7}
	For any natural numbers $n,m$, the $\Psi-$polynomial satisfy the following identity
	\begin{equation}
		\label{WW8}
		\boldmath
		\begin{aligned}
			(2a-b)^{\delta(n)\delta(m)} \: \Psi(a,b,n) \Psi(a,b,m) = \Psi(a,b,n+m) + a^{\min\{n, m\}} \Psi(a,b,n-m).   \\
		\end{aligned}
	\end{equation}
\end{theorem}

\section{\textbf{THE EIGHT LEVELS THEOREM AND FORMULAS FOR $\Psi(a,b,n)$}}

Put $\left( \begin{array}{cc} a & b \\ \alpha & \beta  \end{array} \right) = \left( \begin{array}{cc} 1 & 0 \\ 0 & 1\end{array}\right)$, $(\xi,\eta)=(a,-b)$ in the expansion \eqref{Ready66}, we get the following desirable expansion for polynomial sequence $\Psi(a,b,n)$.  

\begin{theorem}
	\label{expansion-1}
	For any natural number $n$, we have the following expansion 
	\begin{equation}
		\label{expansion-11}
		\begin{aligned}
			\sum_{r=0}^{\lfloor{\frac{n}{2}}\rfloor} (-1)^r  \Psi\left( \begin{array}{cc|r} 1 & 0 & n \\ 0 & 1 & r \end{array} \right)  a^{\lfloor{\frac{n}{2}}\rfloor  - r}  b^r =\Psi(a, b, n).
		\end{aligned}
	\end{equation}
\end{theorem}

For our next step, we will substitute $(a,b)=(xy,-x^2-y^2)$ into Equation \eqref{expansion-11} as per Theorem \eqref{expansion-1}. With the aid of Theorem \eqref{WW3}, we promptly obtain the following outcome:

\begin{theorem}{}
	\label{Expansion-3}  
	For any complex numbers $x,y$, any non negative integers $n,k$, the coefficients $\Psi_k(n)$ of the expansion  
	\begin{equation}
		\label{Expansion-4} 
		\frac{x^n+y^n}{(x+y)^{\delta(n)}} =
		\sum_{k=0}^{\lfloor{\frac{n}{2}}\rfloor}\Psi\left( \begin{array}{cc|r} 1 & 0 & n \\ 0 & 1 & k \end{array} \right) 
		(xy)^{\lfloor{\frac{n}{2}}\rfloor -k} (x^2+y^{2})^{k}
	\end{equation}

\end{theorem}

Depending on $n\pmod 8 $, the Eight Levels theorem, see \cite{2}, gives 8 formulas for the polynomial sequence $\Psi\left( \begin{array}{cc|r} 1 & 0 & n \\ 0 & 1 & r \end{array} \right) $. Now combining the the Eight Levels theorem and Theorem \eqref{expansion-1}, we get the following polynomial expansion for $\Psi(a,b,n)$.

\begin{theorem}{(Expansion for $\Psi(a,b,n)$ as sums of product of difference of squares)}
	\label{Expansion-5}  
	For any complex numbers $x,y$, any non negative integers $n,k$, the coefficients $\Psi_k(n):=  \Psi\left( \begin{array}{cc|r} 1 & 0 & n \\ 0 & 1 & k \end{array} \right)$ of the expansion  
\begin{equation}
	\label{expansion-2}
	\begin{aligned}
	\Psi(a, b, n) =	\sum_{k=0}^{\lfloor{\frac{n}{2}}\rfloor} (-1)^k  \Psi_k(n)  a^{\lfloor{\frac{n}{2}}\rfloor  - k}  b^k.
	\end{aligned}
\end{equation}
	are integers and 
	
	\begin{equation}
		\label{starting} 
		\Psi_0(n) = 		
		\begin{cases}
			+2  &  n \equiv \:\: 0  \pmod{8}  \\   
			+	 1 &  n \equiv \pm 1  \pmod{8}  \\
			\: \: 0  &  n \equiv \pm 2  \pmod{8}  \\   
			-1&  n \equiv \pm 3  \pmod{8}  \\
			-2  &  n \equiv \pm 4  \pmod{8}  \\   
		\end{cases} 
	\end{equation}
	and, for each $ 1 \le k \le  \lfloor{\frac{n}{2}}\rfloor$, the coefficients satisfy the following statements 
	\begin{itemize}	
		\item  For $n \equiv 0, \: 2, \: 4, \: 6\pmod{8}$:\\ 
	\end{itemize}
	
	\underline{$n \equiv 0 \pmod{8}$} 
	\begin{equation*}
		\Psi_k(n) =
		\begin{cases}
			0  & \mbox{for $k$ odd } \\   
			2 \: (-1)^{\lfloor{\frac{k}{2}}\rfloor}
			\: \frac {\prod\limits_{\lambda = 0}^{\lfloor{\frac{k}{2}}\rfloor -1} [n^2 \: - \: (4 \lambda)^2 ]}	{4^k \: k!} & \mbox{for $k$ even } \end{cases}  				
	\end{equation*}
	
	\underline{$n \equiv 2 \pmod{8} $} \\
	\begin{equation*}
		\Psi_k(n) =
		\begin{cases}
			0  & \mbox{for $k$ even } \\   
			2 \: (-1)^{{\lfloor{\frac{k}{2}}\rfloor}} \:
			n \:	\frac {\prod\limits_{\lambda = 1}^{\lfloor{\frac{k}{2}}\rfloor} [ n^2 \: - \: (4 \lambda -2)^2 ]}	{4^k \: k!} & \mbox{for $k$ odd } \end{cases}   
	\end{equation*}
	
	\underline{$n \equiv 4 \pmod{8}$} \\	
	\begin{equation*}	
		\Psi_k(n) =
		\begin{cases}
			0  & \mbox{for $k$ odd } \\   
			2 \: (-1)^{\lfloor{\frac{k}{2}}\rfloor +1}  \:
			\frac {\prod\limits_{\lambda = 0}^{\lfloor{\frac{k}{2}}\rfloor -1} [n^2 \: - \: (4 \lambda)^2 ]}	{4^k \: k!} & \mbox{for $k$ even } \end{cases}   
	\end{equation*}	
	
	\underline{$n \equiv 6 \pmod{8}$}  \\
	\begin{equation*}
		\Psi_k(n) =
		\begin{cases}
			0  & \mbox{for $k$ even } \\   
			2 \: (-1)^{\lfloor{\frac{k}{2}}\rfloor +1}
			\:  n \: \frac {\prod\limits_{\lambda = 1}^{\lfloor{\frac{k}{2}}\rfloor } [n^2 \: - \: (4 \lambda -2)^2 ]} 	{4^k \: k!} & \mbox{for $k$ odd } \end{cases}   
	\end{equation*}
	
	\begin{itemize}	
		\item  For $n \equiv 1, \: 3, \: 5, \: 7\pmod{8}$:\\ 
	\end{itemize}
	
	\underline{$n \equiv 1 \pmod{8}$} 	
	\begin{equation*}
		\Psi_k(n) = (-1)^{\lfloor{\frac{k}{2}}\rfloor} (n+1-2k)^{\delta(k)}  \:
		\frac {\prod\limits_{\lambda = 1}^{\lfloor{\frac{k}{2}}\rfloor} [ (n+1)^2 \: - \: (4 \lambda -2)^2 ]}	{4^k \: k!}  
	\end{equation*}
	
	\underline{$n \equiv 3 \pmod{8}$}  \\		
	\begin{equation*}
		\Psi_k(n) = (-1)^{\lfloor{\frac{k}{2}}\rfloor + \delta(k-1)} (n+1) (n+1-2k)^{\delta(k-1)}  
		\frac {\prod\limits_{\lambda = 1}^{\lfloor{\frac{k-1}{2}}\rfloor} [(n+1)^2 - (4 \lambda)^2 ]}	{4^k \: k!}  
	\end{equation*}		
	
	\underline{$n \equiv 5 \pmod{8}$ }  \\
	\begin{equation*}
		\Psi_k(n) = (-1)^{\lfloor{\frac{k}{2}}\rfloor +1}  \: (n+1-2k)^{\delta(k)} \:
		\frac {\prod\limits_{\lambda = 1}^{\lfloor{\frac{k}{2}}\rfloor }[(n+1)^2 \: - \: (4 \lambda -2)^2 ]}	{4^k \: k!}  
	\end{equation*}
	
	\underline{$n \equiv 7 \pmod{8}$} \\
	\begin{equation*}
		\Psi_k(n) = (-1)^{\lfloor{\frac{k}{2}}\rfloor + \delta(k)} (n+1) (n+1-2k)^{\delta(k-1)}  \:
		\frac {\prod\limits_{\lambda = 1}^{\lfloor{\frac{k-1}{2}}\rfloor } [(n+1)^2 - (4 \lambda)^2 ]}	{4^k \: k!}  \\ 
	\end{equation*}
	
\end{theorem}

\section{\textbf{RESULTS FOR MERSENNE PRIMES AND MERSENNE COMPOSITES}}
Mersenne numbers $2^{p} - 1$ with prime $p$ form the sequence 
\[3, 7, 31, 127, 2047, 8191, 131071, 524287, 8388607, 536870911, \dotsc \] 
(sequence {A001348} in \cite{Slo}). For $2^{p} - 1$ to be prime, it is necessary that $p$ itself be prime.	

\subsection{PRIMALITY TEST FOR MERSENNE PRIMES}
\begin{theorem}{(A new version for Lucas-Lehmer primality test)}
	\label{U14}
	Given prime $p \geq 5$. The number $2^p-1$ is prime if and only if 	
	\begin{equation}
		\label{U15} 
		2n-1 \quad  \vert \quad \Psi(1,4,n),
	\end{equation}
	where $n:=2^{p-1}$.
\end{theorem}

\begin{proof}
	Given prime $p \geq 5$, let $n:=2^{p-1}$. From Lucas-Lehmer test, \cite{Jean}, we have 
	\begin{equation*}
		\begin{aligned}
			2^p -1 \quad \text{is prime} \quad  \iff 2^p -1  \quad | \quad (1+\sqrt{3})^n + (1-\sqrt{3})^n.	
		\end{aligned}
	\end{equation*}
	As $n$ even, $\delta(n)=0$, and from Theorem \eqref{WW3}, we get the following equivalent statement:
	\begin{equation*}
		\begin{aligned}
			2^p -1 \quad \text{is prime} \quad  \iff 2^p -1  \quad | \quad 
			\Psi(x_0 \: y_0,-x_0^2-y_0^2, n),
		\end{aligned}
	\end{equation*}
	where $x_0= 1+\sqrt{3}, \quad  y_0= 1-\sqrt{3}$. As $(x_0 \: y_0,-x_0^2-y_0^2) = (-2, -8),$ and from Theorem \eqref{W11}, and noting $(2^p -1, 2 ) =1 $,  we get the following equivalent statements:
	
	\begin{equation*}
		\begin{aligned}
			2^p -1 \quad \text{is prime} \quad  &\iff 2^p -1  \quad | \quad \Psi(-2,-8,n) \\ 
			&\iff 2^p -1  \quad | \quad(-2)^{\lfloor{\frac{n}{2}}\rfloor}\: \Psi(1,4,n) \\ 	
			&\iff 2^p -1  \quad | \quad  \Psi(1,4,n).  
		\end{aligned}
	\end{equation*}	
\end{proof}

Theorem \eqref{U14} tells us that any Mersenne prime must be a factor for $\Psi(1,4,n)$. Thus we need to study the general arithmetical and combinatorial properties of the sequence $\Psi(1,4,n)$. We need the following theorem.

\begin{theorem}{}
	\label{X-1}
		Given prime $p \geq 5$, and $n:=2^{p-1}$.  $2^p-1$ is prime if and only if for any natural number $\mu$ we have
 \begin{equation}
 	\label{X-2}
 	\Psi(1,4,n\mu) \equiv  
 	\begin{cases}
 		+2 \pmod{2n-1} &  \mu \equiv 0    \pmod{4} \\
  	\; 0 \; \pmod{2n-1} & \mu \equiv 1,3    \pmod{4} \\
  	-2 \pmod{2n-1} &   \mu \equiv 2    \pmod{4} \\
 	\end{cases}.
 \end{equation}
 \end{theorem}
 
 \begin{proof}
 	Let $2^p-1$ be a prime, and $n:=2^{p-1}$. Then from Theorem \eqref{U14} we get $\Psi(1,4,n)  \equiv 0    \pmod{2n-1}$. From Theorem \eqref{WW7}, the following recurrence relation is true for any natural number $\mu$: 
 	
 	\begin{equation}
 		\label{X-3}
 		\begin{aligned}
 			\Psi(1,4,n) \Psi(1,4,n\mu) = \Psi(1,4,n(\mu+1)) + \Psi(1,4,n(\mu-1)).   \\
 		\end{aligned}
 	\end{equation}
Noting that $\Psi(1,4,0)=2$, and from \eqref{X-3}, and working recursively, and starting from $\mu=1$, and noting $\Psi(1,4,n)  \equiv 0    \pmod{2n-1}$, we get the proof of \eqref{X-2} for any natural number $\mu$. Now suppose \eqref{X-2} holds for any for any natural number $\mu$. Then it holds for $\mu =1$, which means $\Psi(1,4,n) \equiv 0    \pmod{2n-1}$, which implies, by Theorem \eqref{U14}, that $2^p-1$ is prime. This completes the proof.
 	
 \end{proof}
  
 	Based on Theorem \eqref{X-1}, and Theorem \eqref{Expansion-5}, we establish the following new result:
  
 \begin{theorem}{(Enhanced Lucas-Lehmer Primality Test)}
 	\label{GG111}
 	For a given prime $p \geq 5$ and $n:=2^{p-1}$, $2^p-1$ is a Mersenne prime if and only if, for any natural number $\mu$, the following congruence holds:
 	\begin{equation}
 		\label{XX-2}
 		\sum_{\substack{k=0,\\ k \:even}}^{\lfloor{\frac{n\mu}{2}}\rfloor}   \quad  (-1)^{\lfloor{\frac{k}{2}}\rfloor}	\: \frac {\prod\limits_{\lambda = 0}^{\lfloor{\frac{k}{2}}\rfloor -1} [(n\mu)^2 - (4 \lambda)^2]}	{ k! } \equiv  
 		\begin{cases}
 			+1 \pmod{2n-1} & \text{if} \quad  \mu \equiv 0 \pmod{4} \\
 			\;\; 0 \; \pmod{2n-1} & \text{if} \quad \mu \equiv 1,3 \pmod{4} \\
 			-1 \pmod{2n-1} & \text{if} \quad \mu \equiv 2 \pmod{4} \\
 		\end{cases}.
 	\end{equation}
 \end{theorem}

 Substituting $\mu=2$ into Equation \eqref{XX-2}, we readily obtain the following necessary condition applicable to all Mersenne primes.
 
 \begin{theorem}{(Necessary Condition for Mersenne Primes)}
 	\label{GGG111}
 	For a given prime $p \geq 5$ and $n:=2^{p-1}$, if $2^p-1$ is a Mersenne prime, then the following congruence must hold:
 	
 	\begin{equation}
 		\label{XXC-2}
 		\sum_{\substack{k=0}}^{\lfloor{\frac{n}{2}}\rfloor}   \quad 	\: \frac {\prod\limits_{\lambda = 0}^{k -1} [ (4 \lambda)^2 - 1]}	{ 2k! } \equiv  
 		-1 \pmod{2n-1}.
 	\end{equation}
 \end{theorem}

\subsection{MERSENNE COMPOSITE NUMBERS	}
The number $2^p-1$ is called Mersenne composite number if $p$ is prime but  $2^p-1$ is not prime.
\begin{theorem}{}
	\label{Maresenn composite}
	Given prime $p$, $n:=2^{p-1}$. If 
	\begin{equation}
		\label{composite} 
		2n-1 \quad  \vert \quad \Psi(1,4,n\: \pm 1),
	\end{equation}
	then $2^p-1$ is a Mersenne composite number.
\end{theorem}

\begin{proof}
	Suppose for contradiction that $2^p - 1$ divides $\Psi(1,4,n \pm 1)$. Now, assume the contrary, that $2^p - 1$ is prime. Then, referring to Theorem \eqref{U14} and using the recurrence relation \eqref{def1} recursively, we arrive at $2^p - 1$ dividing $\Psi(1,4,1) = 1$. This leads to a contradiction.
\end{proof}

\section{\textbf{NEW COMBINATORIAL IDENTITIES RELATED TO MERSENNE NUMBERS}}

From Theorem \eqref{Expansion-5}, we get

\begin{theorem}
		\label{G-B} 
	For any natural number $\tau \equiv 0 \pmod 8$, we get 
	\begin{equation}
		\label{G-B0} 
		2 \; \sum_{k=0}^{\lfloor{\frac{\tau}{4}}\rfloor}  \; \frac {\prod\limits_{\lambda = 0}^{k -1}  [ (4 \lambda)^2 - \tau^2 ]}	{ (2k)! \: 4^{2k}} \; a^{\lfloor{\frac{\tau}{2}}\rfloor  - 2k} \; b^{2k} = \Psi(a,\pm b,\tau). 
	\end{equation}
\end{theorem} 

An immediate result is the following. 
\begin{theorem}
	For any natural number $\tau \equiv 0 \pmod 8$, we get 
	\begin{equation}
		\label{G-I} 
	2 \;	\sum_{k=0}^{\lfloor{\frac{\tau}{4}}\rfloor}    \;   \frac {\prod\limits_{\lambda = 0}^{k -1}  [ (4 \lambda)^2 - \tau^2 ]}	{ (2k)! }   = \Psi(1,\pm 4,\tau). 
	\end{equation}
\end{theorem}

As Mersenne primes factorize the sequence $\Psi(1,4,\tau)$ for specific natural numbers $\tau$, it becomes imperative to thoroughly investigate the associated expansions of $\Psi(1,4,\tau)$. This investigation should extend beyond isolated examinations of $\Psi(1,4,\tau)$; instead, it should encompass a broader context. To attain a comprehensive understanding, it is essential to explore notable expansions of $\Psi(a,b,\tau)$ across all natural numbers $\tau=2^l$, where $l\geq 3$. For example, when considering any natural numbers $\tau=2^l$, $l\geq 3$:
 
\begin{itemize}
	\item Put $(a,b)= (1,1)$ in Theorem \eqref{G-B}:\\
	\begin{equation}
		\label{G-B1} 
		\begin{aligned}				
			2 \; \sum_{k=0}^{\lfloor{\frac{\tau}{4}}\rfloor}  \; \frac {\prod\limits_{\lambda = 0}^{k -1}  [ (4 \lambda)^2 - \tau^2 ]}	{ (2k)! \: 4^{2k}} \; = -1. 
		\end{aligned}
	\end{equation}

	\item Put $(a,b)= (1,2)$ in Theorem \eqref{G-B}:\\
\begin{equation}
	\label{G-B2} 
	\begin{aligned}					
	 \; \sum_{k=0}^{\lfloor{\frac{\tau}{4}}\rfloor}  \; \frac {\prod\limits_{\lambda = 0}^{k -1}  [ (4 \lambda)^2 - \tau^2 ]}	{ (2k)! \: 4^{k}} \;  = 1.	
	\end{aligned}
\end{equation}

\item 

Substitute $(a,b) = (1,\sqrt{2})$ into Theorem \eqref{G-B}:

\begin{equation}
	\label{G-B2} 
	\begin{aligned}					
		\; \sum_{k=0}^{\lfloor{\frac{\tau}{4}}\rfloor}  \; \frac {\prod\limits_{\lambda = 0}^{k -1}  [ (4 \lambda)^2 - \tau^2 ]}	{ (2k)! \: 2^{3k}} \;  = \pm 1,
	\end{aligned}
\end{equation}

where the right-hand side of Equation \eqref{G-B2} equals $+1$ if $\tau \equiv 0 \pmod {16}$ and $-1$ if $\tau \equiv 8 \pmod {16}$.

\end{itemize}

 \section{\textbf{NEW DIFFERENTIAL OPERATOR APPROACH TO MERSENNE PRIMES}}
 
 As we embark on the enigmatic journey of Mersenne primes, a realm of mathematical wonders unfolds before us. Embracing a new differential approach and exploring essential operator lists, we illuminate the intricate properties of Mersenne primes, unraveling their mysteries. The intriguing link between Mersenne primes and other sequences, such as the captivating Lucas sequence, holds the potential for groundbreaking discoveries and novel research avenues. By integrating advanced mathematical techniques like number theory algorithms and computational methods, we deepen our understanding of these elusive numbers and their profound significance. This exploration of Mersenne primes and their interplay with various sequences beckons further investigation and inspires curiosity. The pursuit of knowledge in this domain promises profound insights into the fundamental nature of primes and offers a gateway to innovative advancements in mathematics and beyond. We mustn't overlook the fact that Mersenne primes are nestled within the factors of $\Psi(1,4,n)$. This prompts us to substitute $r=\lfloor{\frac{n}{2}}\rfloor$ into Equation \eqref{diff3} of Theorem \eqref{A1}, in conjunction with Theorem \eqref{exp1}, leading to the following immediate consequence.
 
 \begin{theorem}{(A New Differential Operator Approach to Mersenne prime)}\\
 	\label{GG} For any numbers $a,b$, $4 a -  b \neq 0$, and any natural number $n$, we have
 	\begin{align}
 		\begin{aligned}
 			\frac{1}{(\lfloor{\frac{n}{2}}\rfloor)!} \Big( \frac{{\partial} }{\partial a} + 4 \frac{{\partial}}{\partial b}\Big)^{\lfloor{\frac{n}{2}}\rfloor} \Psi(a,b,n) = \Psi(1,4,n).  
 		\end{aligned}
 	\end{align}
 \end{theorem}

 \begin{theorem}
 	For any natural number $n$, we get:
 	\begin{align}
 		\label{G3} 
 		\frac{1}{(\lfloor{\frac{n}{2}}\rfloor)!} \Big( \frac{{\partial} }{\partial a} \Big)^{\lfloor{\frac{n}{2}}\rfloor} \Psi(a,b,n) &= \Psi(1,0,n),\\
 		\frac{1}{(\lfloor{\frac{n}{2}}\rfloor)!} \Big(  \frac{{\partial}}{\partial b}\Big)^{\lfloor{\frac{n}{2}}\rfloor} \Psi(a,b,n) &= \Psi(0,1,n),\\
 		\frac{1}{(\lfloor{\frac{n}{2}}\rfloor)!} \Big( \frac{{\partial} }{\partial a} + 4 \frac{{\partial}}{\partial b}\Big)^{\lfloor{\frac{n}{2}}\rfloor} \Psi(a,b,n) &= \Psi(1,4,n),\\
 			\frac{1}{(\lfloor{\frac{n}{2}}\rfloor)!} \Big( \frac{{\partial} }{\partial a} + 3 \frac{{\partial}}{\partial b}\Big)^{\lfloor{\frac{n}{2}}\rfloor} \Psi(a,b,n) &= \Psi(1,3,n),\\
 				\frac{1}{(\lfloor{\frac{n}{2}}\rfloor)!} \Big( \frac{{\partial} }{\partial a} +  \frac{{\partial}}{\partial b}\Big)^{\lfloor{\frac{n}{2}}\rfloor} \Psi(a,b,n) &= \Psi(1,1,n),\\
 			\frac{1}{(\lfloor{\frac{n}{2}}\rfloor)!} \Big( \frac{{\partial} }{\partial a} +  2 \frac{{\partial}}{\partial b}\Big)^{\lfloor{\frac{n}{2}}\rfloor} \Psi(a,b,n) &= \Psi(1,2,n),\\ 
 				\frac{1}{(\lfloor{\frac{n}{2}}\rfloor)!} \Big( \frac{{\partial} }{\partial a} + \sqrt{2} \frac{{\partial}}{\partial b}\Big)^{\lfloor{\frac{n}{2}}\rfloor} \Psi(a,b,n) &= \Psi(1,\sqrt{2},n),\\
 			\frac{1}{(\lfloor{\frac{n}{2}}\rfloor)!} \Big( 2 \frac{{\partial} }{\partial a} + 5 \frac{{\partial}}{\partial b}\Big)^{\lfloor{\frac{n}{2}}\rfloor} \Psi(a,b,n) &= \Psi(2,5,n).\\			 			
 	\end{align}
 \end{theorem}

 The theorem ahead establishes bridges between sequences, imparting knowledge about their arithmetic connections.
 
 \begin{theorem}
 		For any $\mu_i, a,b,\alpha_i,\beta_i, n$, $\beta_i a - \alpha_i b \neq 0,$  for each $i$, we get:
 	\begin{align}
 		\label{GQ} 
  \Big[ \sum_{i}^{} \mu_i \Big( \alpha_i \frac{{\partial} }{\partial a} + \beta_i \frac{{\partial}}{\partial b}\Big)^{\lfloor{\frac{n}{2}}\rfloor} \Big ] \Psi(a,b,n) =  (\lfloor{\frac{n}{2}}\rfloor) ! \sum_{i}^{} \mu_i  \Psi(\alpha_i,\beta_i,n).
 	\end{align}
\end{theorem}
 
For instance, this theorem may contribute to our comprehension of the arithmetic interplay between Mersenne primes and the Lucas sequence, and vice versa. If, for example, we were to pinpoint an elegant formula for the sum of the operators
 
 \[      \mu_1 \Big(  \frac{{\partial} }{\partial a} + 4 \frac{{\partial}}{\partial b}\Big)^{m}  + \mu_2  \Big(  \frac{{\partial} }{\partial a} + 3 \frac{{\partial}}{\partial b}\Big)^{m}  + \mu_3 \Big(  \frac{{\partial} }{\partial b} \Big)^{m},    \]

 where $\mu_1, \mu_2, \mu_3$ are integers, and considering the fundamental facts that $\Psi(1,3,m)= \pm L(m)$ and $\Psi(0,1,m)=\pm 1$, we would undoubtedly obtain a new, elegant formula for:
 
 \[
 \mu_1\Psi(1,4,m) \pm \mu_2 L(m) \pm \mu_3.
 \]

\section{\textbf{UNIFICATION OF WELL-KNOWN SEQUENCES}}\label{Unify}

The polynomials $\Psi(a,b,n)$ exhibit captivating arithmetical and differential properties, uniting a plethora of well-known polynomials and sequences. Notably, these include Mersenne numbers, the Chebyshev polynomials of the first and second kind, Dickson polynomials of the first and second kind, Lucas numbers, Fibonacci numbers, Fermat numbers, Pell-Lucas polynomials, Pell numbers, and others (see \cite{222}, \cite{333}). These remarkable polynomial sequences form a rich tapestry that connects diverse areas of mathematics, providing valuable insights and deepening our understanding of these well-established sequences.

\subsection{\textbf{CONNECTIONS WITH FIBONACCI AND LUCAS SEQUENCES}}

Notably, the $\Psi$ and $\Phi$ polynomials are closely linked to the Fibonacci and Lucas sequences through a variety of formulas, which are presented below:

\begin{equation*}
	\begin{aligned}
		\Psi(-1,-3,n) &= L(n), \qquad \Psi(1,-3,n) &=
		\begin{cases}
			F(n) & \text{for $n$ odd}, \\
			L(n) & \text{for $n$ even}. \end{cases} \\
	\end{aligned}
\end{equation*}

Here, $L(n)$ represents the Lucas numbers defined by the recurrence relation $L(0)=2$, $L(1)=1$, and $L(n+1)=L(n)+L(n-1)$, while $F(n)$ denotes the Fibonacci numbers defined by the recurrence relation $F(0)=0$, $F(1)=1$, and $F(n+1)=F(n)+F(n-1)$. These formulas showcase the remarkable connections between the $\Psi$ and the well-known Fibonacci and Lucas sequences. Moreover, for any natural number $n$, we get, from \eqref{def1}, the following relation

\begin{equation}
	\label{(1,Root5)}
	\Psi(1,\sqrt{5},n) = 
	\begin{cases}
		L(\frac{n}{2}) &   n \equiv 0    \pmod{4} \\
		L(\frac{n+1}{2}) + F(\frac{n-1}{2})\sqrt{5} &  n \equiv 1 \pmod{4} \\
		- F(\frac{n}{2}) \sqrt{5}  &    n \equiv 2 \pmod{4} \\
		-L(\frac{n-1}{2})-F(\frac{n+1}{2})\sqrt{5}& n \equiv 3\pmod{4} 
	\end{cases}.
\end{equation}

and

\begin{equation}
	\label{(1,Root5)}
	\Psi(1,-\sqrt{5},n) = 
	\begin{cases}
		L(\frac{n}{2}) &   n \equiv 0    \pmod{4} \\
		L(\frac{n+1}{2}) - F(\frac{n-1}{2})\sqrt{5} &  n \equiv 1 \pmod{4} \\
		 F(\frac{n}{2}) \sqrt{5}  &    n \equiv 2 \pmod{4} \\
		-L(\frac{n-1}{2})+F(\frac{n+1}{2})\sqrt{5}& n \equiv 3\pmod{4} 
	\end{cases}.
\end{equation}

\subsection{\textbf{CONNECTIONS WITH MERSENNE AND FERMAT NUMBERS}}

Notably, the choices $(2,-5)$ and $(-2,-5)$ for $(a,b)$ establish significant links between the $\Psi$-Polynomial and Mersenne numbers. For each natural number $n$, the following desirable formulas emerge:

\begin{equation}
	\label{(2,-5)-(-2,-5)}
	\begin{aligned}
		\Psi(-2,-5,n) &= 2^{n} +(-1)^{n} &,& \qquad \Psi(2,-5,n) &= \frac{2^{n} +1}{3^{\delta(n)}}. 
	\end{aligned}
\end{equation}

For $p$ odd, we obtain $ \Psi(-2,-5,p) = 2^{p} -1 = M_{p}$, where $M_{p}$ are known as Mersenne numbers. Mersenne primes $M_p$, for some prime $p$, are remarkable due to their connection with primality tests and perfect numbers. In number theory, a perfect number is a positive integer equal to the sum of its proper positive divisors, excluding the number itself. It remains unknown whether any odd perfect numbers exist. Recent research \cite{25} has shown that odd perfect numbers, if they exist, must be greater than $10^{1500}$. Additionally, we observe that $\Psi(-2,-5,2^n) = \Psi(2,-5,2^n) = 2^{2^n} + 1 = F_{n}$, where $F_{n}$ represents the Fermat number. The connections with Mersenne and Fermat numbers further emphasize the significance and broad applications of the $\Psi$ polynomials in the realm of number theory.

\subsection{\textbf{CONNECTIONS WITH PELL NUMBERS AND POLYNOMIALS}}

The Pell numbers $P_n$ are defined by the recurrence relation $P_0 = 0$, $P_1 = 1$, and $P_n = 2P_{n-1} + P_{n-2}$. 

The Pell-Lucas numbers $Q_n$ are defined by the recurrence relation $Q_0 = 2$, $Q_1 = 2$, and $Q_n = 2Q_{n-1} + Q_{n-2}$. We find that the $\Psi$ polynomial serves as a generalization for the Pell-Lucas numbers: $Q_{n} = 2^{\delta(n)} \Psi(-1,-6,n)$. Similarly, the Pell-Lucas numbers $Q_n = Q_n(1)$, where $Q_n(x)$ is the Pell-Lucas polynomial defined by the recurrence relation $Q_0(x) = 2$, $Q_1(x) = 2x$, and $Q_{n+1}(x) = 2xQ_{n}(x) + Q_{n-1}(x)$. This leads to the following relation:

\begin{equation*}
	Q_{n}(x) = (2x)^{\delta(n)} \; \Psi(-1,-2-4x^2,n).
\end{equation*}

\subsection{\textbf{CONNECTIONS WITH DICKSON POLYNOMIALS}}

The Dickson polynomial of the first kind of degree $n$ with parameter $\alpha$, denoted $D_n(x,\alpha)$, see \cite{222, 333, L-5, 24B, L-6}, are defined by the following formulas, where $\delta(n)=1$ for $n$ odd and $\delta(n)=0$ for $n$ even:

\begin{equation*}
	\begin{aligned}
		D_n(x,\alpha) &= \sum_{i=0}^{\left\lfloor \frac{n}{2} \right\rfloor}\frac{n}{n-i} \binom{n-i}{i} (-\alpha)^i x^{n-2i}.
	\end{aligned}
\end{equation*}

Considering Equation \eqref{comp3} and Theorem \eqref{comp2}, it becomes evident that the following relation holds between the $\Psi$ polynomials and the Dickson polynomials:

\begin{equation*}
	\begin{aligned}
		D_n(x,\alpha) &= x^{\delta(n)} \Psi(\alpha,2\alpha-x^2,n).
	\end{aligned}
\end{equation*}

\subsection{\textbf{CONNECTIONS WITH CHEBYSHEV POLYNOMIALS}}

The Chebyshev polynomial of the first kind of degree $n$, denoted $T_n(x)$, see \cite{Chebyshev, 6, 222, 333, L-5, Cody, 22, 8,  9}, are defined by the following formulas:

\begin{equation*}
	\begin{aligned}
		T_n(x) &= \sum_{i=0}^{\left\lfloor \frac{n}{2} \right\rfloor}(-1)^i \frac{n}{n-i} \binom{n-i}{i} (2)^{n-2i-1} x^{n-2i}.
	\end{aligned}
\end{equation*}

In light of Equation \eqref{comp3} and Theorem \eqref{comp2}, it is manifest that the following connection is established between the $\Psi$ polynomials and the Chebyshev polynomials:

\begin{equation*}
	\begin{aligned}
		T_n(x) &= \frac{x^{\delta(n)}}{2^{\delta(n+1)}} \Psi(1,2-4x^2,n).
	\end{aligned}
\end{equation*}

\subsubsection{$\boldsymbol{\Psi(1,2,n)}$ - \textbf{ Unraveling Captivating Patterns with $a=1, b=2$}}

In the realm of $\Psi$-sequences, a mesmerizing revelation awaits as we delve into the enchanting world of $\Psi(1,2,n)$ with $a=1, b=2$. This captivating case leads us to a striking formula:

\[ \Psi(1,2,n)=  \:  \frac{n}{n-\left\lfloor \frac{n}{2} \right\rfloor} \binom{n-\left\lfloor \frac{n}{2} \right\rfloor}{\left\lfloor \frac{n}{2} \right\rfloor} \:(-1)^{\left\lfloor \frac{n}{2} \right\rfloor}.      \]

As we explore further, the beauty of this sequence unfolds, revealing a profound expression that showcases the symmetries and connections hidden within:

\[ \Psi(1,2,n)= \:(-1)^{\left\lfloor \frac{n}{2} \right\rfloor} \:  2^{\delta(n-1)} \:n^{\delta(n)}.      \]

The allure of this case lies in its ability to uncover the intricate interplay of numbers, offering us glimpses into the underlying structures and relationships within the sequence.

\section{\textbf{SYMMETRICAL BEAUTY: PERIODIC PATTERNS IN $\boldsymbol{\Psi}$ SEQUENCES}}
\subsection{Revealing Additional Special Cases of the $\boldsymbol{\Psi}$-Sequence}

Delving deeper into the world of the $\Psi$-sequence, we unearth a wealth of other exceptional scenarios where distinct values of $a$ and $b$ intertwine to form this fascinating class of polynomials. Each special case holds its own unique set of patterns and properties, enriching our understanding of this newly discovered mathematical structure.

\subsubsection{$\boldsymbol{\Psi(1,1,n)}$}

For $a=1$ and $b=1$, the $\Psi$-sequence takes the following form:
\begin{equation}
	\label{PP0}
	\Psi(1,1,n) =
	\begin{cases}
		+2 & n \equiv \pm 0 \pmod{6} \\
		+1 & n \equiv \pm 1 \pmod{6} \\
		-1 & n \equiv \pm 2 \pmod{6} \\
		-2 & n \equiv \pm 3 \pmod{6}
	\end{cases}
\end{equation}

This specific case of the $\Psi$-sequence provides valuable insights into the patterns and characteristics of the polynomial for the given values of $a$ and $b$.

\subsubsection{$\boldsymbol{\Psi(1,0,n)}$ - Unraveling Mersenne Prime Patterns}

The fascinating world of the $\Psi$-sequence offers yet another captivating revelation when $a=1$ and $b=0$. In this special case, the $\Psi$-sequence takes the form:
\begin{equation}
	\label{PP0Q}
	\Psi(1,0,n) =
	\begin{cases}
		+2 & n \equiv \pm 0 \pmod{8} \\
		+1 & n \equiv \pm 1 \pmod{8} \\
		\: 0 & n \equiv \pm 2 \pmod{8} \\
		-1 & n \equiv \pm 3 \pmod{8} \\
		-2 & n \equiv \pm 4 \pmod{8} \\
	\end{cases}
\end{equation}

This alluring special case unravels intriguing patterns related to Mersenne primes, offering a glimpse into the enigmatic world of these unique numbers.

\subsubsection{$\boldsymbol{\Psi(1,-1,n)}$ - Unraveling Patterns with $a=1, b=-1$}

Yet another compelling case within the realm of the $\Psi$-sequence emerges when $a=1$ and $b=-1$. In this instance, the $\Psi$-sequence takes on the following form:

\begin{equation}
	\label{peroidicity1}
	\Psi(1,-1,n) =
	\begin{cases}
		+2 & n \equiv \pm 0 \pmod{12} \\
		+1 & n \equiv \pm 1 , \pm 2 \pmod{12} \\
		\:0 & n \equiv \pm 3 \pmod{12} \\
		-1 & n \equiv \pm 4, \pm 5 \pmod{12} \\
		-2 & n \equiv \pm 6 \pmod{12}
	\end{cases}
\end{equation}

This specific case provides us with invaluable insights into the periodic patterns of the $\Psi$-sequence when $a=1$ and $b=-1$. The graph of $\Psi(1,-1,n)$ exhibits a mesmerizing symmetry, evoking a delicate beauty reminiscent of a flower. With a periodicity of 12, the pattern gracefully repeats, creating an enchanting interplay of numbers. This captivating visualization provides a delightful glimpse into the underlying nature of the sequence.

\begin{figure}[ht]
	\centering
	\begin{tikzpicture}[->,>=stealth',shorten >=0.1pt,auto,node distance=2.5cm,
		thick,main node/.style={circle,draw,font=\sffamily\Large\bfseries,
			minimum size=3mm}]
		\node[main node] (A) {\bf \huge 2};
		\node[main node] (B) [left of=A] {1};
		\node[main node] (C) [left of=B] {0};
		\node[main node] (D) [left of=C] {-1};
		\node[main node] (E) [left of=D] {-2};
		\path[every node/.style={font=\sffamily\small,
			fill=white,inner sep=1pt}]
		(A) edge [bend left=30] node[below=1mm] {\bf \Large 1} (B)
		(B) edge [out=120,in=60,out distance= 3 cm,in distance=4cm] node[above=1mm] {2} (B)
		(B) edge [bend left=30] node[below=1mm]{3} (C)
		(C) edge [bend left=30] node[below=1mm]{4} (D)
		(D) edge [out=120,in=60,out distance= 3 cm,in distance=4cm] node[above=1mm] {5} (D)
		(D) edge [bend left=30] node[below=1mm]{6} (E)
		(E) edge [bend left=30] node[above=1mm] {7} (D)
		(D) edge [out=300,in=240,out distance= 3 cm,in distance=4cm] node[below=1mm] {8} (D)
		(D) edge [bend left=30] node[above=1mm] {9} (C)
		(C) edge [bend left=30] node[above=1mm] {10} (B)
		(B) edge [out=300,in=240,out distance= 3 cm,in distance=4cm] node[below=1mm] {11} (B)
		(B) edge [bend left=30] node[above=1mm] {12} (A);
	\end{tikzpicture}
	\caption{Symmetric and Elegant Graph of $\Psi(1,-1,n)$ with Periodicity = 12}
	\label{Rfig2}
\end{figure}
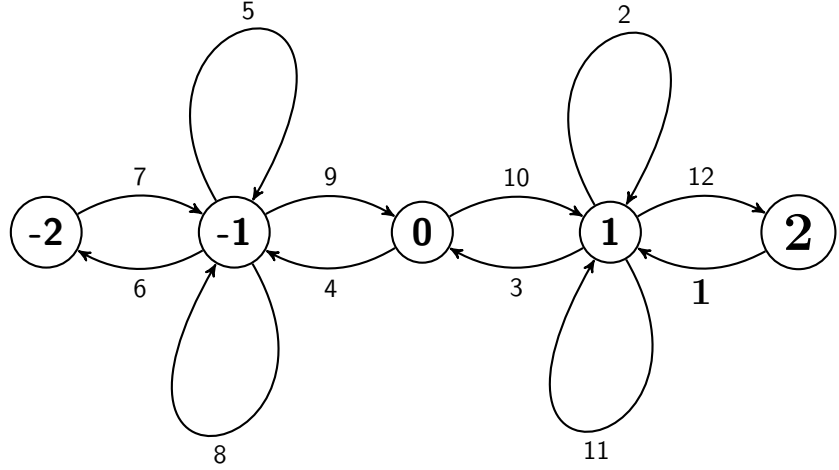

\subsubsection{$\boldsymbol{\Psi(1,\sqrt{2},n)}$ }
It is quite unexpected to note that the periodicity of $\Psi(1,\sqrt{2},n)$ equals $16$ and its formula is given by:

\begin{equation}
	\label{root-2}
	\Psi(1,\sqrt{2} ,n) =
	\begin{cases}
		2 & n \equiv \pm 0 \pmod{16} \\
		1 & n \equiv \pm 1 \pmod{16} \\
		-\sqrt{2} & n \equiv \pm 2 \pmod{16} \\
		-1-\sqrt{2} & n \equiv \pm 3, \pmod{16} \\
		0 & n \equiv \pm 4 \pmod{16} \\
		1+\sqrt{2} & n \equiv \pm 5, \pmod{16} \\
		\sqrt{2} & n \equiv \pm 6 \pmod{16} \\
		-1 & n \equiv \pm 7 \pmod{16} \\
		-2 & n \equiv \pm 8 \pmod{16} \\
	\end{cases}
\end{equation}

The graph of $\Psi(1,\sqrt{2},n)$ exhibits mesmerizing symmetry with a periodicity of 16, forming patterns that strikingly resemble parts of a flower. Each element in the sequence is connected in a harmonious manner, reminiscent of the blooming flowers in a vast garden. The intriguing properties of this graph and its connection to the golden ratio encourage further investigation into the captivating behavior of $\Psi$ polynomials and their links to fundamental mathematical constants. The graph of $\Psi(1,\sqrt{2},n)$ is as follows:

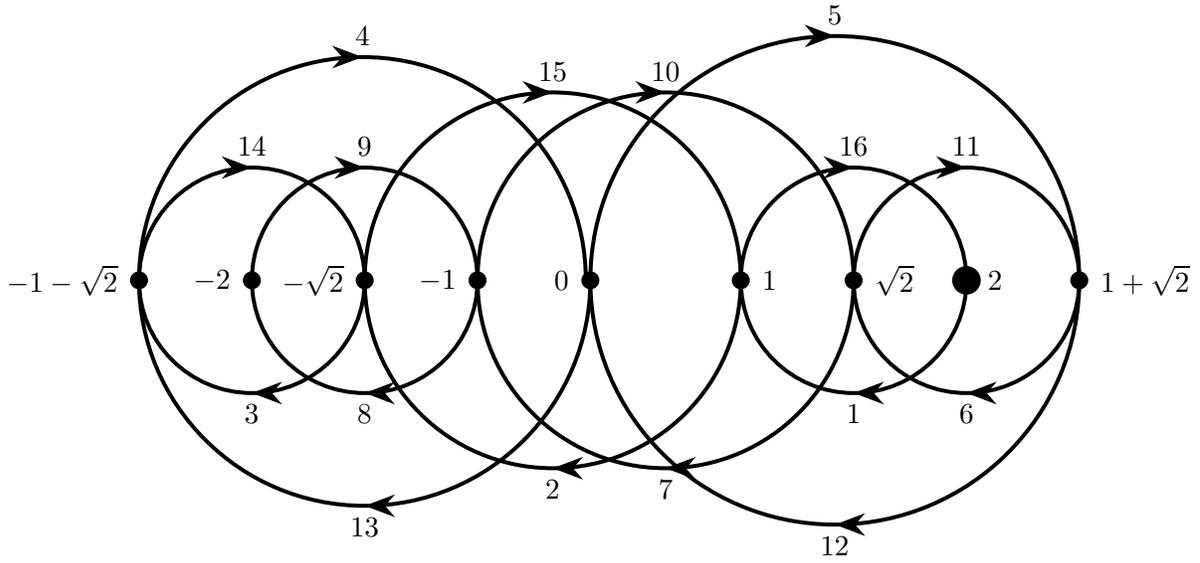
\begin{figure}[ht]
	\centering
	\begin{tikzpicture}[>={Stealth[length=4mm]}, node distance=1.8cm, thick, main node/.style={circle, draw, font=\sffamily\small\bfseries, minimum size=1cm},
		halfcirc/.style={postaction={decorate,decoration={markings,mark=at position 0.5 with {\arrow{>}}}}}]
		% Nodes on the left
		\foreach \x/\label in {-8/{$-1-\sqrt{2}$}, -6.5/{$-2$}, -5/{$-\sqrt{2}$}, -3.5/{$-1$}, -2/{$0$}}
		\filldraw (\x,0) circle (3pt) node[label=left:{\label}] {};
		
		% Nodes on the right
		\foreach \x/\label in {4.5/{$1+\sqrt{2}$}, 3/{$2$}, 1.5/{$\sqrt{2}$}, 0/{$1$}}
		\filldraw (\x,0) circle (3pt) node[label=right:{\label}] {};
		
		% Half circles as edges
		\draw[halfcirc, line width=1.5pt] (3,0) arc[start angle=360, end angle=180, radius=1.5cm] node[midway, below] {1};
		\draw[halfcirc, line width=1.5pt] (0,0) arc[start angle=360, end angle=180, radius=2.5cm] node[midway, below] {2};
		\draw[halfcirc, line width=1.5pt] (-5,0) arc[start angle=360, end angle=180, radius=1.5cm] node[midway, below] {3};
		
		\draw[halfcirc, line width=1.5pt] (-8,0) arc[start angle=180, end angle=0, radius=2.97cm] node[midway, above] {4};
		
		\draw[halfcirc, line width=1.5pt] (-2,0) arc[start angle=180, end angle=0, radius=3.25cm] node[midway, above] {5};
		
		% Arrow and label for the last edge
		\draw[halfcirc, line width=1.5pt] (4.5,0) arc[start angle=360, end angle=180, radius=1.5cm] node[midway, below] {6};

		% Arrow and label for the last edge
		\draw[halfcirc, line width=1.5pt] (1.5,0) arc[start angle=360, end angle=180, radius=2.5cm] node[midway, below] {7};

		% Arrow and label for the last edge
		\draw[halfcirc, line width=1.5pt] (-3.5,0) arc[start angle=360, end angle=180, radius=1.5cm] node[midway, below] {8};

		\draw[halfcirc, line width=1.5pt] (-6.5,0) arc[start angle=180, end angle=0, radius=1.5cm] node[midway, above] {9};

		\draw[halfcirc, line width=1.5pt] (-3.5,0) arc[start angle=180, end angle=0, radius=2.5cm] node[midway, above] {10};
		
		\draw[halfcirc, line width=1.5pt] (1.5,0) arc[start angle=180, end angle=0, radius=1.5cm] node[midway, above] {11};

		% Arrow and label for the last edge
		\draw[halfcirc, line width=1.5pt] (4.5,0) arc[start angle=360, end angle=180, radius=3.25cm] node[midway, below] {12};
		
		% Arrow and label for the last edge
		\draw[halfcirc, line width=1.5pt] (-2,0) arc[start angle=360, end angle=180, radius=3cm] node[midway, below] {13};

		\draw[halfcirc, line width=1.5pt] (-8,0) arc[start angle=180, end angle=0, radius=1.5cm] node[midway, above] {14};

		\draw[halfcirc, line width=1.5pt] (-5,0) arc[start angle=180, end angle=0, radius=2.5cm] node[midway, above] {15};

		\draw[halfcirc, line width=1.5pt] (0,0) arc[start angle=180, end angle=0, radius=1.5cm] node[midway, above] {16};
		
		% Special style for node "2"
		\filldraw (3,0) circle (5pt);
		
	\end{tikzpicture}
	\caption{The graph of $\Psi(1,\sqrt{2},n)$, periodicity = 16}
	\label{fig-1}
\end{figure}

\newpage
\subsubsection{$\boldsymbol{\Psi(1,\sqrt{3},n)}$}
The periodicity of $\Psi(1,\sqrt{3},n)$ equals $24$;

\begin{equation}
	\label{root-3}
	\Psi(1,\sqrt{3} ,n) = 
	\begin{cases}
		2	      &   n \equiv \pm 0             \pmod{24} \\
		1	      &   n \equiv \pm 1,\pm 4              \pmod{24} \\
		-\sqrt{3} &   n \equiv \pm 2             \pmod{24} \\
		-1-\sqrt{3}	  &   n \equiv \pm 3,        \pmod{24} \\
		2+\sqrt{3}       &   n \equiv \pm 5             \pmod{24} \\
		0             &   n \equiv \pm 6,       \pmod{24} \\
		-2-\sqrt{3}  &   n \equiv \pm 7             \pmod{24} \\
		-1	      &   n \equiv \pm 8, \pm 11             \pmod{24} \\
		1+\sqrt{3}	      &   n \equiv \pm 9            \pmod{24} \\
		\sqrt{3}      &   n \equiv \pm 10,       \pmod{24} \\
		-2  &   n \equiv \pm 12             \pmod{24} \\
	\end{cases}
\end{equation}

\newpage

\begin{figure}[ht]
	\centering
	\begin{tikzpicture}[>={Stealth[length=4mm]}, thick, main node/.style={circle, draw, font=\sffamily\small\bfseries, minimum size=1cm},
		halfcirc/.style={postaction={decorate,decoration={markings,mark=at position 0.5 with {\arrowreversed{>}}}}}]

		% Nodes on the left
		\foreach \y/\label in {9.5/{$2+\sqrt{3}$}, 8/{$1+\sqrt{3}$}, 6.5/{$2$}, 5/{$\sqrt{3}$}, 3.5/{$1$}, 2/{$0$}}
		\filldraw (0,\y) circle (3pt) node[above right] {\label};

		\foreach \y/\label in {0/{$-1$}, -1.5/{$-\sqrt{3}$}, -3/{$-2$}, -4.5/{$-1-\sqrt{3}$},-6/{$-2-\sqrt{3}$ }}
		\filldraw (0,\y) circle (3pt) node[below right] {\label};

		% Half circles as edges
		\draw[halfcirc, line width=1.5pt] (0,3.5) arc[start angle=270, end angle=90, radius=1.5cm] node[midway, left] {1};

			\draw[halfcirc, line width=1.5pt] (0,-1.5) arc[start angle=270, end angle=90, radius=2.5cm] node[pos=0.58, left] {2};

		\draw[halfcirc, line width=1.5pt] (0,-4.5) arc[start angle=270, end angle=90, radius=1.5cm] node[midway, left] {3};
		
		\draw[halfcirc, line width=1.5pt] (0,3.5) arc[start angle=90, end angle=-90, radius=4 cm] node[midway, right] {4};
		
		\draw[halfcirc, line width=1.5pt] (0,9.5) arc[start angle=90, end angle=-90, radius=3 cm] node[pos=0.56, right] {5};
		
		\draw[halfcirc, line width=1.5pt] (0,2) arc[start angle=270, end angle=90, radius=3.75cm] node[midway, left] {6};
		
		\draw[halfcirc, line width=1.5pt] (0,-6) arc[start angle=270, end angle=90, radius=4cm] node[midway, left] {7};
		
		\draw[halfcirc, line width=1.5pt] (0,0) arc[start angle=90, end angle=-90, radius=3cm] node[midway, right] {8};
		
		\draw[halfcirc, line width=1.5pt] (0,8) arc[start angle=90, end angle=-90, radius=4cm] node[midway, right] {9};
		
		\draw[halfcirc, line width=1.5pt] (0,5) arc[start angle=270, end angle=90, radius=1.5cm] node[midway, left] {10};
		
		\draw[halfcirc, line width=1.5pt] (0,0) arc[start angle=270, end angle=90, radius=2.5cm] node[pos=0.54, left] {11};
		
		\draw[halfcirc, line width=1.5pt] (0,-3) arc[start angle=270, end angle=90, radius=1.5cm] node[midway, left] {12};
		
		\draw[halfcirc, line width=1.5pt] (0,0) arc[start angle=90, end angle=-90, radius=1.5cm] node[midway, right] {13};
		
		\draw[halfcirc, line width=1.5pt] (0,5) arc[start angle=90, end angle=-90, radius=2.5cm] node[midway, right] {14};
		
		\draw[halfcirc, line width=1.5pt] (0,8) arc[start angle=90, end angle=-90, radius=1.5cm] node[midway, right] {15};
		
		\draw[halfcirc, line width=1.5pt] (0,0) arc[start angle=270, end angle=90, radius=4cm] node[midway, left] {16};
		
		\draw[halfcirc, line width=1.5pt] (0,-6) arc[start angle=270, end angle=90, radius=3cm] node[pos=0.6, left] {17};
		
		\draw[halfcirc, line width=1.5pt] (0,2) arc[start angle=90, end angle=-90, radius=4cm] node[midway, right] {18};
		
		\draw[halfcirc, line width=1.5pt] (0,9.5) arc[start angle=90, end angle=-90, radius=3.75cm] node[pos=0.53, right] {19};
		
		\draw[halfcirc, line width=1.5pt] (0,3.5) arc[start angle=270, end angle=90, radius=3cm] node[pos=0.55, left] {20};
		
		\draw[halfcirc, line width=1.5pt] (0,-4.5) arc[start angle=270, end angle=90, radius=4cm] node[midway, left] {21};
		
		\draw[halfcirc, line width=1.5pt] (0,-1.5) arc[start angle=90, end angle=-90, radius=1.5cm] node[midway, right] {22};

		\draw[halfcirc, line width=1.5pt] (0,3.5) arc[start angle=90, end angle=-90, radius=2.5cm] node[pos=0.6, right] {23};
		
	\draw[halfcirc, line width=1.5pt] (0,6.5) arc[start angle=90, end angle=-90, radius=1.5cm] node[midway, right] {24};
	
		% Special style for node "2"
	\filldraw (0,6.5) circle (5pt);

	\end{tikzpicture}
	\caption{The graph of $\Psi(1,\sqrt{3},n)$, periodicity = 24}
\end{figure}
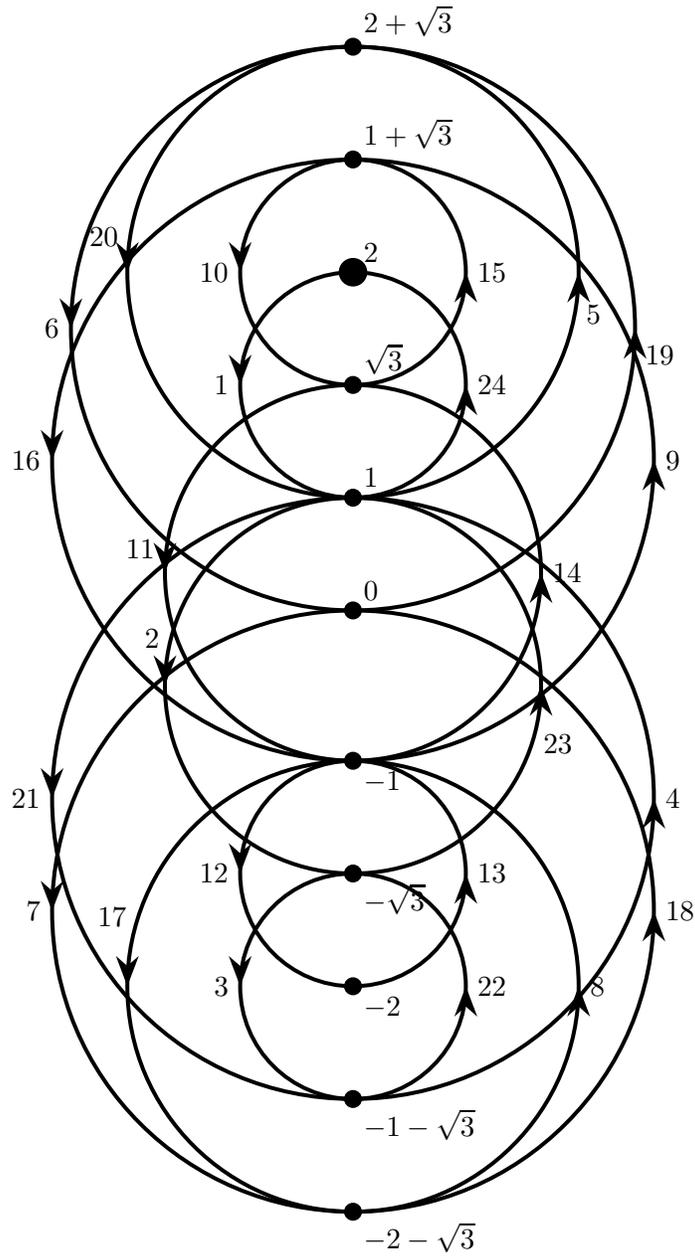

	\subsubsection{\textbf{PERIODICITY OF $\Psi(1,\phi - 1,n)$, LINKS WITH THE GOLDEN RATIO}}
	
	Notably, the $\Psi$ Polynomials reveal intriguing links with the Golden Ratio. Specifically, when $a=1$ and $b=\phi-1$, we obtain the following expressions for the $\Psi$ and $\Phi$ sequences:
	
	\begin{equation}
		\label{Golden ratio1}
		\Psi(1,\phi - 1,n) =
		\begin{cases}
			2 & n \equiv \pm 0 \pmod{20} \\
			1 & n \equiv \pm 1 \pmod{20} \\
			-\phi +1 & n \equiv \pm 2 \pmod{20} \\
			
			-\phi & n \equiv \pm 3, \pm 4 \pmod{20} \\
			0 & n \equiv \pm 5 \pmod{20} \\
			\phi & n \equiv \pm 6, \pm 7 \pmod{20} \\
			\phi - 1 & n \equiv \pm 8 \pmod{20} \\
			-1 & n \equiv \pm 9 \pmod{20} \\
			-2 & n \equiv \pm 10 \pmod{20} \\
		\end{cases}
	\end{equation}
	
	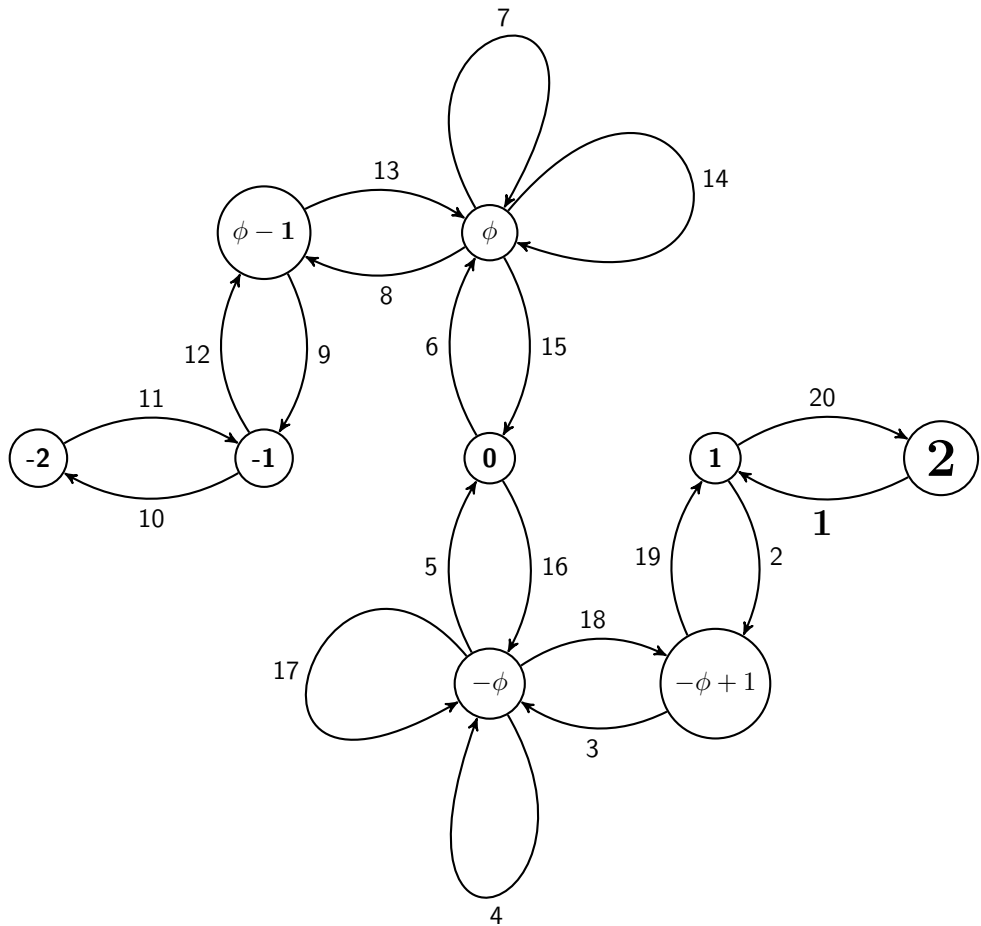
\begin{figure}[ht]
		\centering
		\begin{tikzpicture}[->,>=stealth',shorten >=0.1pt,auto,auto,node distance=3cm,
			thick,main node/.style={circle,draw,font=\sffamily\small\bfseries,
				minimum size=0.1mm}]
			\node[main node] (A) {\bf \huge 2};
			\node[main node] (B) [left of=A] {1};
			\node[main node] (C) [left of=B] {0};
			\node[main node] (D) [left of=C] {-1};
			\node[main node] (E) [left of=D] {-2};
			\node[main node] (M) [above of=C] {$\phi$};
			\node[main node] (N) [left of=M] {$\bf \phi-1$};
			\node[main node] (U) [below of=C] {$-\phi$};
			\node[main node] (V) [right of=U] {$-\phi+1$};
			\path[every node/.style={font=\sffamily\small,
				fill=white,inner sep=1pt}]
			% Right-hand-side arrows rendered from top to bottom to
			% achieve proper rendering of labels over arrows.
			(A) edge [bend left=30] node[below=1mm] {\bf \Large{1}} (B)
			(B) edge [bend left=30] node[right=1mm]{2} (V)
			(V) edge [bend left=30] node[below=1mm]{3} (U)
			(U) edge [out=300,in=250,out distance= 3 cm,in distance=4cm] node[below=1mm] {4} (U)
			(U) edge [bend left=30] node[left=1mm] {5} (C)
			(C) edge [bend left=30] node[left=1mm] {6} (M)
			(M) edge [out=120,in=60,out distance= 3 cm,in distance=4cm] node[above=1mm] {7} (M)
			(M) edge [bend left=30] node[below=1mm] {8} (N)
			(N) edge [bend left=30] node[right=1mm] {9} (D)
			(D) edge [bend left=30] node[below=1mm] {10} (E)
			(E) edge [bend left=30] node[above=1mm] {11} (D)
			(D) edge [bend left=30] node[left=1mm] {12} (N)
			(N) edge [bend left=30] node[above=1mm] {13} (M)
			(M) edge [out=50,in=340,out distance= 4 cm,in distance=4cm] node[right=1mm] {14} (M)
			(M) edge [bend left=30] node[right=1mm] {15} (C)
			(C) edge [bend left=30] node[right=1mm] {16} (U)
			(U) edge [out=130,in=210,out distance= 3 cm,in distance=4cm] node[left=1mm] {17} (U)
			(U) edge [bend left=30] node[above=1mm] {18} (V)
			(V) edge [bend left=30] node[left=1mm] {19}(B)
			(B) edge [bend left=30] node[above=1mm] {20} (A);
		\end{tikzpicture}
		\caption{Graph of $\Psi(1,\phi-1,n)$ with periodicity = 20}
		\label{fig1}
	\end{figure}
	
	The graph of $\Psi(1,\phi-1,n)$ gracefully repeats with a periodicity of 20, forming elegant patterns that strikingly resemble parts of a flower. Each element in the sequence is connected in a harmonious manner, much like flowers blooming in a vast garden, inviting further exploration into the captivating behavior of the $\Psi$ sequences. This remarkable symmetry and linkage with the Golden Ratio beckon researchers to delve deeper into the enchanting world of these polynomials and their connections to fundamental mathematical constants.

\section{\textbf{A WINDOW INTO THE FUTURE OF MERSENNE PRIMES}}
\subsection{\textbf{BEYOND THE SURFACE: UNEXPLORED HORIZONS IN MERSENNE PRIMES AND POLYNOMIAL SEQUENCES}}

The domain of Mersenne prime research is both profound and teeming with promise, harboring the potential for countless future discoveries. As we conclude our paper, my intention is to provide a concise glimpse of a pivotal result, tailored to ignite the curiosity and enthusiasm of our readers for the ever-evolving realm of Mersenne primes. Our final revelation unveils a remarkable new finding within the Mersenne prime landscape, with exciting implications that beckon further exploration. We introduce the double-indexed polynomial sequences $A_r(k)$ and $B_r(k)$ and divulge their intriguing ratios, as illustrated below:
\begin{equation}
	\label{UU}
	\begin{aligned}
		\frac{A_0( \lfloor{\frac{p}{2}}\rfloor )}{(p-1)(p-2) \cdots (p - \lfloor{\frac{p}{2}}\rfloor )} \quad , \quad \frac{B_0( \lfloor{\frac{n}{2}}\rfloor )}{(n-1)(n-2) \cdots (n - \lfloor{\frac{n}{2}}\rfloor )}.
	\end{aligned}
\end{equation}

The intriguing revelation is that these ratios yield integers.

Remarkably, we uncover a direct correlation between the primality of $2^p - 1$ and the divisibility relationships among these ratios. This intriguing association calls for deeper exploration, promising valuable insights into the distribution and essence of Mersenne primes. Our findings set the stage for promising advancements in the field of number theory and ignite the spark for future inquiries into the underlying mathematical structure. Succinctly, without delving into the proof, we present the following striking new result for Mersenne primes, ushering in exciting possibilities for further investigation.

\begin{theorem}{}
	\label{Theorem of ABC1}
	For any given prime $p\geq 5$, let $n:=2^{p-1}$, and we associate the double-indexed polynomial sequences $A_r(k)$ and $B_r(k)$, defined as follows:
	\begin{equation}
		\label{ABC2} 
		\begin{aligned}		
			A_r(k) &=  (p-r-k) \: A_r(k-1) + \: 4 \: (p-2r) \: A_{r+1}(k-1), \quad &A_r(0) = 1  \quad   \text{for all} \:\: r,\\
			B_r(k) &=  -2 \: (n-r-k) \: B_r(k-1) - \:2 \: (n-2r-1) \: B_{r+1}(k-1),  \quad &B_r(0) = 1  \quad   \text{for all} \:\: r.
		\end{aligned}
	\end{equation}
	Then both of the ratios   
	\begin{equation}
		\label{ABC3} 
		\begin{aligned}
			\frac{A_0( \lfloor{\frac{p}{2}}\rfloor   )}{(p-1)(p-2) \cdots (p - \lfloor{\frac{p}{2}}\rfloor )} \quad	, \quad 	\frac{B_0( \lfloor{\frac{n}{2}}\rfloor   )}{(n-1)(n-2) \cdots (n - \lfloor{\frac{n}{2}}\rfloor )} 	
		\end{aligned}
	\end{equation}
	are integers. Moreover the number $2^p -1$ is prime \bf{if and only if} 
	\begin{equation}
		\label{ABC4} 
		\begin{aligned}
			\frac{A_0( \lfloor{\frac{p}{2}}\rfloor   )}{(p-1)(p-2) \cdots (p - \lfloor{\frac{p}{2}}\rfloor )}  \quad  \vert \quad \frac{B_0( \lfloor{\frac{n}{2}}\rfloor   )}{(n-1)(n-2) \cdots (n - \lfloor{\frac{n}{2}}\rfloor )}.
		\end{aligned}
	\end{equation}
	
\end{theorem}

\subsection{\textbf{CONCLUSION}}
In conclusion, this paper has presented a new result for Mersenne primes, shedding light on intriguing connections between the double-indexed polynomial sequences $A_r(k)$ and $B_r(k)$. The ratios described in Equation \eqref{UU} were shown to yield integers, and a significant link between the primality of $2^p - 1$ and the divisibility relation of these ratios was established. With rigorous investigation and additional work, researchers can explore the properties and behaviors of the polynomial sequences $A_r(k)$ and $B_r(k)$, potentially unveiling deeper insights into the distribution and nature of Mersenne primes. The ``if and only if''  statement in Theorem (\ref{Theorem of ABC1}), linking the divisibility of the ratios, presents an intriguing direction for future investigations. The existence of integer ratios of \eqref{ABC3} hints at a hidden mathematical structure yet to be fully understood. Further investigation into the properties of these polynomial sequences and their relationship to prime numbers could lead to significant advancements in the field of number theory and unveil new insights into the distribution and existence of Mersenne primes.

\end{document}